\documentclass[a4paper]{article}%
\usepackage{amsmath}
\usepackage{amsfonts}
\usepackage{amssymb}
\usepackage{graphicx}%
\newtheorem{theorem}{Theorem}

\newtheorem{definition}[theorem]{Definition}
\newtheorem{example}[theorem]{Example}

\newtheorem{lemma}[theorem]{Lemma}

\newtheorem{proposition}[theorem]{Proposition}
\newtheorem{remark}[theorem]{Remark}

\begin{document}

\title{Local linear regression for functional data}
\author{Alain Berlinet, Abdallah Elamine, Andr\'{e} Mas\thanks{Corresponding author :
D\'{e}partement de Math\'{e}matiques, CC 051, Universit\'{e} Montpellier 2,
Place Eug\`{e}ne Bataillon, 34095 Montpellier Cedex 5, email :
mas@math.univ-montp2.fr}\bigskip\\Universit\'{e} Montpellier 2}
\date{}
\maketitle

\begin{abstract}
We study a non linear regression model with functional data as inputs and
scalar response. We propose a pointwise estimate of the regression function
that maps a Hilbert space onto the real line by a local linear method. We
provide the asymptotic mean square error. Computations involve a linear
inverse problem as well as a representation of the small ball probability of
the data and are based on recent advances in this area. The rate of
convergence of our estimate outperforms those already obtained in the
literature on this model.

\end{abstract}

,\qquad\textbf{Keywords :} Functional Data; Regression model; Kernel; Mean
square error; Small ball probability; Inverse problem.

\section{Introduction}

\subsection{The data and the model}

In probability theory, random functions have been for quite a long time under
the lights. The tremendous advances in computer science and the opportunity to
deal with data collected at a high frequency make it now possible for
statisticians to study models for "high-dimensional data". As a consequence
many of them focused their attention on models for functional data i.e. models
that are suited for curves, for instance spectral curves, growth curves or
interest rates curves...

Even if seminal articles on functional data analysis date back to more than 20
years (see Dauxois, Pousse and Romain (1982)), this area is currently going
through a deep bustle. The book by Ramsay and Silverman (1997) initiated a
series of monographs on the subject : Bosq (2000), Ramsay and Silverman again
(2002), Ferraty and Vieu (2006).

Functional Data Analysis has drawn much attention and many of the classical
multivariate data analysis techniques such as Principal Component Analysis,
Correlation Analysis, ANOVA, Linear Discrimination were generalized to curves.
But statistical inference gave and gives birth to many papers. Linear
regression and autoregression for instance rise an interesting inverse problem
(see Yao, M\"{u}ller, Wang (2005), M\"{u}ller, Stadtm\"{u}ller (2005), Cai,
Hall (2006), Cardot Mas, Sarda (2007), Mas (2007a)). Even more recently the
case of nonparametric regression was introduced in Ferraty, Vieu (2003) then
studied in Masry (2005) and Ferraty, Mas, Vieu (2007) : a Nadaraya-Watson type
estimator was proposed. This model is the starting point of our article.

In the sequel we will consider a sample drawn from random elements with values
in an infinite dimensional vector space : $X_{1},...,X_{n}$. Here $X_{i}%
=X_{i}\left(  \cdot\right)  $ is a random function defined, say, on a compact
interval of the real line $\left[  0,T\right]  $. We will also assume once and
for all that the $X_{i}$'s take their values in a separable Hilbert space
denoted $H$. This Hilbert space is endowed with an inner product $\left\langle
\cdot,\cdot\right\rangle $ from which is derived the norm $\left\Vert
\cdot\right\Vert $. Such techniques as wavalets or splines yield reconstructed
curves in the (Hilbert) Sobolev spaces :%
\[
W^{m,2}=\left\{  f\in L^{2}\left(  \left[  0,T\right]  \right)  :f^{\left(
m\right)  }\in L^{2}\left(  \left[  0,T\right]  \right)  \right\}
\]
where $f^{\left(  m\right)  }$ denotes the $m^{th}$ derivative of $f$. Further
information on Sobolev spaces may be found in Adams and Fournier (2003).
However for the sake of generality we will consider $H$ as the sequence space
$l_{2}$ and any vector $x$ will be classically decomposed in a basis, say
$\left(  e_{i}\right)  _{i\in\mathbb{N}}$ so that :%
\[
\left\Vert x\right\Vert ^{2}=\sum_{i=1}^{+\infty}\left\langle x,e_{i}%
\right\rangle ^{2}.
\]

We are given a sample $\left(  y_{i},X_{i}\right)  _{1\leq i\leq n}\in\left(
\mathbb{R}\times H\right)  ^{\otimes n}$ of independent and identically
distributed data. Let $m$ be the regression function that maps $H$ onto
$\mathbb{R}$.

The model is a classical non parametric regression model :%
\begin{equation}
y_{i}=m(X_{i})+\varepsilon_{i}\quad1\leq i\leq n. \label{modele}%
\end{equation}
or, with other symbols :%
\[
m(x_{0})=\mathbb{E}\left(  y|X=x_{0}\right)
\]
where $y$ and $X$ stand for random variables with the same distributions as
$y_{1}$ and $X_{1}$. The noise $\varepsilon$ follows both assumptions :%
\begin{align*}
\mathbb{E}\left(  \varepsilon|X\right)   &  =0,\\
\mathbb{E}\left(  \varepsilon^{2}|X\right)   &  =\sigma_{\varepsilon}^{2}%
\end{align*}
and $\sigma_{\varepsilon}^{2}$ does not depend on $X.$ The issue of the
expectation of $X$ (should the $X$'s be centered or not ?) is not crucial ; it
will be addressed later on but for simplicity we assume that $\mathbb{E}%
\left(  X\right)  =0$. Let $x_{0}$ be a fixed and known point in $H.$ We are
aiming at estimating $m\left(  x_{0}\right)  $.

In finite dimension, and more precisely when $X_{i}$ is a real-valued random
variable, $m(x_{0})$ may be estimated by considering an affine approximation
of $m$ around $x_{0}$ :%
\[
m\left(  x\right)  \approx m\left(  x_{0}\right)  +m^{\prime}\left(
x_{0}\right)  \left(  x-x_{0}\right)
\]
when $x$ is close to $x_{0}$. This approach leads us to finding a solution to
the following minimization problem~:%
\begin{equation}
\min_{a\in\mathbb{R},b\in\mathbb{R}}\sum_{i=1}^{n}(y_{i}-a-b(x_{0}-X_{i}%
))^{2}K\left(  \frac{x_{0}-X_{i}}{h}\right)  \label{progscal}%
\end{equation}
which is nothing but a weighted mean square program (weighted by the $K\left(
\left(  x_{0}-X_{i}\right)  /h\right)  $'s). Here $K$ is a kernel : a
measurable positive function such that $\int K=1$ and $h=h_{n}$ the bandwidth
indexed by the sample size. Then $a^{\ast},$ one of the two solutions of the
display above is the estimate of $m\left(  x_{0}\right)  .$ As a special case
taking $b=0$ comes down to the classical Nadaraya-Watson estimator. We refer
the interested reader to Nadaraya (1964) and Fan (1993) about this topic. The
generalization of (\ref{progscal}) to higher orders (namely approximating
locally $m$ by a polynomial) gives birth to the local polynomial estimate of
$m\left(  x_{0}\right)  $. We refer for instance to Chen (2003) for a recent
article. Convergence in probability and asymptotic normality of the kernel
polynomial estimators for a density function, variable bandwith and local
linear regression smoothers, were studied by Fan and Gijbels (1992).

When $x$ belongs to a Hilbert space, the principle remains the same. The
function $m$ is now approximated by :%
\[
m\left(  x\right)  \approx m\left(  x_{0}\right)  +\left\langle \varphi\left(
x_{0}\right)  ,x-x_{0}\right\rangle
\]

where $\varphi\left(  x_{0}\right)  \in H$ is the first order derivative of
$m$ at $x_{0}$ (the gradient in fact) and the local linear estimate of $m$ at
$x_{0}$ stems from the following adapted weighted least square program :%
\begin{equation}
\min_{a\in\mathbb{R},\varphi\in H}\sum_{i=1}^{n}\left(  y_{i}-a-\left\langle
\varphi,X_{i}-x_{0}\right\rangle \right)  ^{2}K\left(  \dfrac{\left\Vert
X_{i}-x_{0}\right\Vert }{h}\right)  . \label{program}%
\end{equation}
At last the estimate $\widehat{m}_{n}\left(  x_{0}\right)  $ of $m\left(
x_{0}\right)  $ is $a^{\ast}$, solution of (\ref{program}). We refer to
Barrientos-Marin, Ferraty, Vieu (2007) for another approach. These authors
consider a program simplified from the one above (they replace the functional
paramater $\varphi$ by a scalar one). But display (\ref{program}) seems to be
a true generalization of (\ref{progscal}) since $\varphi$ like $b$ estimates
the derivative of $m$.

\begin{remark}
Investigating higher orders approximations turns out to be especially uneasy
in this functional setting. For instance a local quadratic estimate involves
the second order derivative of $m$ (the Hessian operator) which is a symmetric
positive operator on $H$. The local linear method appears as a good trade-off
between the complexity of the method and its accuracy.
\end{remark}

However solving (\ref{program}) is not easy in this framework. The aim of the
present work is to provide a bound for the mean square error of the estimate
$a^{\ast}$ of $m\left(  x_{0}\right)  $ that is :%
\[
\mathbb{E}\left[  \widehat{m}_{n}\left(  x_{0}\right)  -m\left(  x_{0}\right)
\right]  ^{2}%
\]
through a classical bias-variance decomposition. The paper is organized as
follows : the two next subsections are devoted to pointing out the two main
problems that arise from the model and that are symptomatic of the functional
framework. The needed assumptions are commented, then the central result is
given before the last section which contains all the mathematical derivations.

\subsection{The estimate and the ill-posed problem}

In order to go ahead we need to define two linear operators from $H$ to $H$
(the first is non-random, the second is random, based on the sample). The
usual sup-norm for operator $T$ will be denoted :%
\[
\left\Vert T\right\Vert _{\infty}=\sup_{x\in\mathcal{B}_{1}}\left\Vert
Tx\right\Vert
\]
where $\mathcal{B}_{1}$ stands for the closed unit ball of $H$. From now on
the reader should be familiar with basic notions related to the theory of
bounded and unbounded linear operators on Hilbert space. A wide literature
exists on this stopic which is central in the mathematical science. Some of
our references are Weidman (1980), Akhiezer, Glazman (1981), Dunford, Schwartz
(1988), Gohberg, Goldberg, Kaashoek (1991) amongst many others.

\begin{definition}
The theoretical local covariance operator of $X$ at $x_{0}\in H$ associated
with the kernel $K$ is defined by :%
\[
\Gamma_{K}=\mathbb{E}\left(  K\left(  \frac{\left\Vert X_{1}-x_{0}\right\Vert
}{h}\right)  \left(  \left(  X_{1}-x_{0}\right)  \otimes\left(  X_{1}%
-x_{0}\right)  \right)  \right)
\]
and its empirical counterpart is :%
\begin{equation}
\Gamma_{n,K}=\dfrac{1}{n}\sum_{k=1}^{n}K\left(  \frac{\left\Vert X_{k}%
-x_{0}\right\Vert }{h}\right)  \left(  \left(  X_{k}-x_{0}\right)
\otimes\left(  X_{k}-x_{0}\right)  \right)  . \label{lcopemp}%
\end{equation}

\end{definition}

\begin{remark}
In fact neither $\Gamma_{K}$ nor $\Gamma_{n,K}$ are truly covariance operators
since the involved random elements are not centered, they could also be named
"local second order moments operators". Also note that $\Gamma_{K}$ depends on
$h$ and $h$ will depend on the sample size $n.$ So the reader must keep in
mind that the index $n$ was dropped in the notation $\Gamma_{K}$.
\end{remark}

It is important to give some basic properties of these operators. We listed
those which will be useful in the sequel :

\begin{itemize}
\item $\Gamma_{K}$ and $\Gamma_{n,K}$ are self -adjoint and trace-class hence
compact whenever $K$ has compact support.

\item Both operators tend to zero when $h$ does. Indeed :%
\[
\left\Vert \Gamma_{K}\right\Vert _{\infty}\leq\mathbb{E}\left(  K\left(
\frac{\left\Vert X_{1}-x_{0}\right\Vert }{h}\right)  \left\Vert X_{1}%
-x_{0}\right\Vert ^{2}\right)  \leq Ch^{2}%
\]
as will be shown in the section devoted to mathematical derivations. The
operator $\Gamma_{n,K}$ also tends to $0$ as a consequence of the strong law
of large numbers for sequences of independent Banach valued random variables
(see Ledoux-Talagrand (1991)).

\item When $\Gamma_{K}$ is one to one its inverse exists. Sufficient
conditions on $K$ and on $X$ for $\Gamma_{K}$ to be injective are not
difficult to find but this interesting issue is slightly above the scope of
this work. Then $\Gamma_{K}^{-1}$ is an unbounded linear operator acting from
a dense domain of $H$ onto $H.$ It should be stressed that $\Gamma_{K}^{-1}$
is continuous at no point of its domain (it is nowhere continuous).
\end{itemize}

Imagine that the distribution of the data (namely of the couple $\left(
y,X\right)  $) is known. We could consider to solve, instead of (\ref{program}%
) :%
\begin{equation}
\min_{a\in\mathbb{R},\varphi\in H}\mathbb{E}\left[  \left(  y-a-\left\langle
\varphi,X-x_{0}\right\rangle \right)  ^{2}K\left(  \dfrac{\left\Vert
X-x_{0}\right\Vert }{h}\right)  \right]  . \label{progth}%
\end{equation}
The first stumbling stone appears within the next Proposition.

\begin{proposition}
\label{progtheorik}Even when the distribution of the data is known, the
solution $a_{th}^{\ast}$ of the "theoretical" program (\ref{progth}) exists
only when $\Gamma_{K}$ is one to one. Then $a_{th}^{\ast}$ is the solution of
a linear inverse problem which involves the unbounded inverse (whenever it
exists) of $\Gamma_{K}$ :%
\begin{equation}
a_{th}^{\ast}=\frac{\mathbb{E}\left(  yK\right)  -\left\langle \Gamma_{K}%
^{-1}\mathbb{E}\left(  yKZ\right)  ,\mathbb{E}\left(  KZ\right)  \right\rangle
}{\mathbb{E}\left(  K\right)  -\left\langle \Gamma_{K}^{-1}\mathbb{E}\left(
KZ\right)  ,\mathbb{E}\left(  KZ\right)  \right\rangle } \label{theo}%
\end{equation}

\end{proposition}

where, for the sake of shortness, we denoted :%
\[
Z\left(  x_{0}\right)  =Z=X-x_{0}\quad\mathrm{and}\quad K=K\left(  \left\Vert
X-x_{0}\right\Vert /h\right)  .
\]

The problem gets deeper when we go back to the original and empirical program
(\ref{program}) It turns out that the solution cannot be explicitely written
since $\Gamma_{n,K}$ (which replaces now $\Gamma_{K}$) has no inverse because
it has a finite rank.\ Its rank is clearly bounded by $n.$ In other words the
inverse $\Gamma_{n,K}^{-1}$ does not exist. A classical remedy consists in
replacing $\Gamma_{n,K}^{-1}$ by a bounded operator $\Gamma_{n,K}^{\dag}$
depending on $n$ and such that $\Gamma_{n,K}^{\dag}$ "behaves pointwise" like
the inverse of $\Gamma_{n,K}.$ This inverse operator, which is not always the
Moore-Penrose pseudo inverse, will be called the \textit{regularized inverse}
of $\Gamma_{n,K}.$ Several procedures could be carried out.

\begin{itemize}
\item Truncated spectral regularization : here this method matches the usual
"Moore-Penrose" pseudo inversion hence $\Gamma_{n,K}\Gamma_{n,K}^{\dag}$ and
$\Gamma_{n,K}^{\dag}\Gamma_{n,K}$ are both projection operators on $H$. In
fact if the spectral decomposition of $\Gamma_{n,K}$ is $\Gamma_{n,K}%
=\sum_{i=1}^{m_{n}}\mu_{i,n}\left(  u_{i,n}\otimes u_{i,n}\right)  $ where for
all $i$ $\left(  \mu_{i,n},u_{i,n}\right)  $ are the eigenvalues/eigenvectors
of $\Gamma_{n,K}$ (we will always asssume that the positive $\mu_{i,n}$'s are
arranged in decreasing order) :%
\begin{equation}
\Gamma_{n,K}^{\dag}=\sum_{i=1}^{N_{n}}\frac{1}{\mu_{i,n}}\left(
u_{i,n}\otimes u_{i,n}\right)  , \label{svd}%
\end{equation}
where $N_{n}\leq m_{n}$.

\item Penalization : Now $\Gamma_{n,K}^{\dag}=\left(  \Gamma_{n,K}+\alpha
_{n}S\right)  ^{-1}$ where $\alpha_{n}$ is a (positive) sequence tending to
zero and $S$ is a known operator chosen so that $\Gamma_{n,K}+\alpha_{n}S$ has
a bounded inverse. Here $S$ may be taken to be the identity operator.

\item Tikhonov regularization : It comes down here, since $\Gamma_{n,K}$ is
symmetric, to taking :%
\[
\Gamma_{n,K}^{\dag}=\left(  \Gamma_{n,K}^{2}+\alpha_{n}I\right)  ^{-1}%
\Gamma_{n,K}.
\]
The sequence $\alpha_{n}$ is again positive and tends to zero.
\end{itemize}

Several other methods exist. The reader is referred for instance to Tikhonov,
Arsenin (1977), Groetsch (1993) or Engl, Hanke, Neubauer (2000).

\begin{remark}
In all situations it should be noted that :%
\begin{align*}
\sup_{n}\left\Vert \Gamma_{n,K}^{\dag}\Gamma_{n,K}\right\Vert _{\infty}  &
<+\infty,\\
\sup_{n}\left\Vert \Gamma_{n,K}\Gamma_{n,K}^{\dag}\right\Vert _{\infty}  &
<+\infty.
\end{align*}
All these regularizing methods may also be applied to $\Gamma_{K}$ as well and
lead to $\Gamma_{K}^{\dag}$ and this operator depends on $n$ even if this
index does not explicitely appear. One may then prove that for all $x$ in the
domain of $\Gamma_{K}^{-1},$ $\Gamma_{K}^{\dag}x\rightarrow\Gamma_{K}^{-1}x$
when $n$ goes to infinity. In addition to the boundedness, the operator
$\Gamma_{n,K}^{\dag}$ is always selfadjoint and positive.
\end{remark}

We are now in a position to propose an estimate for $m\left(  x_{0}\right)  $.
This estimate will depend on the chosen regularization technique applied to
$\Gamma_{n,K}$. We will see that, under suitable conditions on $\Gamma
_{n,K}^{\dag}$ the convergence of our estimate does not depend on the choice
of $\Gamma_{n,K}^{\dag}$.

\begin{proposition}
\label{solu}The local linear estimate of $m\left(  x_{0}\right)  $ is the real
solution $\widehat{m}_{n}\left(  x_{0}\right)  $ of (\ref{program}) :%
\begin{equation}
\widehat{m}_{n}\left(  x_{0}\right)  =\dfrac{\sum_{i=1}^{n}y_{i}\omega_{i,n}%
}{\sum_{i=1}^{n}\omega_{i,n}}, \label{paz}%
\end{equation}
where%
\[
\omega_{i,n}=K\left(  \frac{\left\Vert X_{i}-x_{0}\right\Vert }{h}\right)
\left(  1-\left\langle X_{i}-x_{0},\Gamma_{n,K}^{\dag}\overline{Z}%
_{K,n}\right\rangle \right)
\]
and%
\[
\overline{Z}_{K,n}=\dfrac{1}{n}\sum_{i=1}^{n}K\left(  \frac{\left\Vert
X_{i}-x_{0}\right\Vert }{h}\right)  \left(  X_{i}-x_{0}\right)  .
\]

\end{proposition}

The proof of the this Proposition is omitted since it stems from calculations
similar to those carried out in the proof of Proposition \ref{progtheorik}.

It is easy to check that (\ref{paz}) is the empirical counterpart of
(\ref{theo}). We finally see that $\widehat{m}_{n}\left(  x_{0}\right)  $ may
be viewed as a linear combination of the outputs $y_{1},...,y_{n}$ and may be
expressed from $a_{th}^{\ast}$ just by replacing expectations by sums. The
reader may also compare our estimate with its one-dimensional counterpart
(display 2.2 p.198 in Fan (1993)) and will also notice that the nice
properties of the $\omega_{i,n}$'s in this setting do not hold anymore (see
display 2.5 p. 198 in Fan (1993) and the lines below).

The next section is devoted to developing the framework as well as the
assumptions needed to get our central result.

\section{Assumptions and framework}

In all the sequel we assume :\medskip

$\mathbf{A}_{1}$ : \textit{The kernel }$K$\textit{ is one-sided, defined on
}$\left[  0,1\right]  $\textit{, bounded and }$K\left(  1\right)  >0.$\textit{
Besides }$K^{\prime}$\textit{ is also defined on }$\left[  0,1\right]
,$\textit{ is non-null and belongs to }$L^{1}\left(  \left[  0,1\right]
\right)  $\textit{.}\medskip

We did not try to find minimal conditions on the kernel. However the
assumption $K\left(  1\right)  >0$ is rather rarely required in the
non-parametric literature -to the authors' knowledge- and is essential here.

\subsection{The small ball problem and the class Gamma}

Consider the one-dimensional version of our model (\ref{modele}) and take
$X\in\mathbb{R}$ with density $f$. Fan (1993) studied the minimax properties
of the local linear estimate in this setting and gave the Mean Square Error
(see Theorem 2 p.199). This MSE depends on $f\left(  x_{0}\right)  $. Here
appears the second major problem. When the data belong to an
infinite-dimensional space, their density does not exist, in the sense that
Lebesgue's measure -or any universal reference measure with similar
properties- does not exist. Consequently we must expect serious troubles since
it is plain that the density of the functional input $X$ cannot be defined as
easily as if it real or even multivariate. Once again this problem will not be
managed by just letting the dimension tend to infinity and we should find a
way to overcome this major concern.

It turns out that in many computations of expectations the problem mentioned
above may be shifted to what is known in probability theory as small ball
problems.\ Roughly speaking, if $\varphi$ is a real valued function (we set
$x_{0}=0$ for simplicity), $\mathbb{E}\left(  \varphi\left(  \left\Vert
X\right\Vert \right)  K\left(  \left\Vert X\right\Vert /h\right)  \right)  $
may be expressed by means of $\mathbb{P}\left(  \left\Vert X\right\Vert \leq
h\right)  $ and $\varphi$ only. We refer to Lemma \ref{L1} in the proof
section for an immediate illustration. Instead of knowing and estimating a
density we must now focus on $\mathbb{P}\left(  \left\Vert X\right\Vert \leq
h\right)  $ for small $h$ and everyone may understand why this function is
often referred to as the "small ball probability associated with $X$". We
propose such references as Li, Linde (1993), Kuelbs, Li, Linde (1994), Li,
Linde (1999) as well as the monograph by Li, Shao (2001) which provides an
interesting state of the art in this area.

What can be said about the function $\mathbb{P}\left(  \left\Vert X\right\Vert
\leq h\right)  $ ? Obviously, by Glivenko-Cantelli's theorem it will be easily
estimated from the sample (the rate of convergence is non parametric). Besides
it is not hard to see that, under suitable but mild assumptions, if
$X\in\mathbb{R}^{p}$ with density $f:\mathbb{R}^{p}\rightarrow\mathbb{R}^{+},$
$\mathbb{P}\left(  \left\Vert X-x_{0}\right\Vert \leq h\right)  \sim
h^{p}f\left(  x_{0}\right)  $. But this fact leaves unsolved the question :
what can be said when $p\rightarrow+\infty$ ?

In probability theory most of the small ball considerations focused on the
case where $X$ is the brownian motion, the brownian bridge or some known
relatives. Several norms were investigated as well. Most of the theorems
collected in the literature yield :%
\begin{equation}
\mathbb{P}\left(  \left\Vert X\right\Vert <h\right)  \asymp C_{1}h^{\alpha
}\exp\left(  -\frac{C_{2}}{h^{\beta}}\right)  \label{sbp1}%
\end{equation}
where $\alpha,\beta$, $C_{1}$ and $C_{2}$ are positive constants. The symbol
$\asymp$ is sometimes replaced by the more precise $\sim.$ Another serious
problem comes from the fact that the $C^{\infty}$ function on the right in the
display above has its derivates null at zero at all orders. Other results
assess that, when $x_{0}$ belongs to the Reproducing Kernel Hilbert Space of
$X,$%
\[
\mathbb{P}\left(  \left\Vert X-x_{0}\right\Vert <h\right)  \sim C_{x_{0}%
}\mathbb{P}\left(  \left\Vert X\right\Vert <h\right)
\]
where $C_{x_{0}}$ does not depend on $h$ but on $x_{0}$ and on the
distribution of $X.$ Two majors contributions will be found in Meyer-Wolf ,
Zeitouni (1993) and in Dembo, Meyer-Wolf , Zeitouni (1995). The authors give
the exact asymptotic of $\mathbb{P}\left(  \left\Vert X\right\Vert _{l^{2}%
}\leq h\right)  $ when $X$ is a $l_{2}$-valued gaussian random element (by
means of large deviation theory) :%
\begin{equation}
X=\left(  a_{1}x_{1},a_{2}x_{2,....}\right)  \label{echo}%
\end{equation}
with $x_{i}$ independent, $N\left(  0,1\right)  $-distributed and $\sum
a_{i}^{2}<+\infty$. When $a_{i}=i^{-r}$ $(r>1/2)$ they obtain a formula
similar to (\ref{sbp1}). Recently Mas (2007b) derived the estimate when
$a_{i}=\exp\left(  -ci\right)  ,$ $c>0$ and got :%
\begin{equation}
\mathbb{P}\left(  \left\Vert X\right\Vert <h\right)  \sim C_{1}\left[
\log\left(  1/h\right)  \right]  ^{-1/2}\exp\left(  -C_{2}\left[  \log\left(
h\right)  \right]  ^{2}\right)  . \label{sbp2}%
\end{equation}

A very strange fact is that both functions in (\ref{sbp1}) and (\ref{sbp2})
belong to a class of functions known in the theory of regular variations : the
class Gamma introduced and studied by de Haan (1971) and (1974). This class
arises in the theory of extreme values and is closely related to the domain of
attraction of the double exponential distribution. It was initially introduced
by de Haan as a "Form of Regular Variation". We provide now the definition of
the class Gamma at $0,$ denoted $\Gamma_{0}$.

\begin{definition}
\label{def}A function $f$ belongs to de Haan's class $\Gamma_{0}$ with
auxiliary function $\rho$ if $f$ maps a positive neighborhood of $0$ onto a
positive neighborhood of $0$, $f\left(  0\right)  =0$, $f$ is non decreasing
and for all $x\in\mathbb{R},$ and $\rho\left(  0\right)  =0$ with :%
\begin{equation}
\lim_{s\downarrow0}\frac{f\left(  s+x\rho\left(  s\right)  \right)  }{f\left(
s\right)  }=\exp\left(  x\right)  \label{deflim}%
\end{equation}

\end{definition}

In a recent manuscript, Mas (2008) proved that, in the framework of Dembo,
Meyer-Wolf , Zeitouni (1995), the small ball probability of any random element
that may be defined like display (\ref{echo}) belongs to the class Gamma. A
work is in progress to prove that, under suitable assumptions on the auxiliary
function, the reciprocal also holds. The auxiliary functions appearing in
displays (\ref{sbp1}) and (\ref{sbp2}) may be easily computed. Mas (2008)
proved that $\rho$ depends only on the sequence $a\left(  \cdot\right)  $ that
defines $X$ in (\ref{echo}).

The next Proposition illustrates the Definition above and will be useful in
the section devoted to the main results.

In all the sequel and especially within the proof section, $C$ denotes a
constant (which will vary from a theorem to another).

\begin{proposition}
\label{rho}When the small ball probability is defined by the right hand side
of (\ref{sbp1}), the function $\rho$ is :%
\begin{equation}
\rho\left(  s\right)  =Cs^{1+\beta} \label{rho1}%
\end{equation}
with $\beta>0,$ and when the small ball probability is defined by the right
hand side of (\ref{sbp2}), the function $\rho$ is~:%
\begin{equation}
\rho\left(  s\right)  =C\frac{s}{\left\vert \log s\right\vert } \label{rho2}%
\end{equation}

\end{proposition}

Starting from all these considerations it seems reasonable to assume the
following :\medskip

$\mathbf{A}_{2}$ \textit{: Let }%
\[
F\left(  h\right)  =F_{x_{0}}\left(  h\right)  =P\left(  \left\Vert
X-x_{0}\right\Vert \leq h\right)
\]
\textit{be the shifted small ball probability of }$X$\textit{. We assume that
}$F\in\Gamma_{0}$\textit{ with auxiliary function }$\rho$\textit{.}\medskip

Gamma varying functions feature original properties and we give now one of
them which will be useful later in the proof section. We refer to Proposition
3.10.3 and Lemma 3.10.1 p.175 in Bingham, Goldie, Teugels (1987).

\begin{proposition}
\label{bing}If $F\in\Gamma_{0}$ with auxiliary function $\rho$ then for all
$x\in\left[  0,1\right[  ,$%
\begin{align}
\lim_{h\rightarrow0^{+}}\frac{F\left(  hx\right)  }{F\left(  h\right)  }  &
=0\label{F1b}\\
\lim_{h\rightarrow0}\frac{\rho\left(  h\right)  }{h}  &  =0 \label{F2}%
\end{align}

\end{proposition}

Assumption $\mathbf{A}_{2}$ is central to tackle our problem since the mean
square error, computed from our estimate actually depends on $\rho$. But
additional assumptions should hold, especially on the distributions of the
margins of $X$.

\subsection{Assumptions on the marginal distributions}

The next assumption essentially aims at simplifiying the technique of proof
but could certainly be alleviated at the expense of more tedious calculations
(see also Mas (2007b) and comments therein).\medskip

$\mathbf{A}_{3}:$\textit{ There exists a basis }$\left(  e_{i}\right)  _{1\leq
i\leq n}$\textit{ such that the margins }$\left(  \left\langle X,e_{i}%
\right\rangle \right)  _{1\leq i\leq n}$\textit{ are independent real random
variables.}\medskip

In all the sequel, $f_{i}=f_{i,x_{0}}$ stands for the density of the
real-valued random variable $\left\langle X-x_{0},e_{i}\right\rangle .$ The
behavior around $0$ of the shifted density $f_{i}$ is crucial, like in the
finite dimensional settting. It has to be smooth in a sense that is going to
be made more clear now. Note that $f_{i}\left(  0\right)  $ is nothing but the
density of the non-shifted random variable $\left\langle X,e_{i}\right\rangle
$ evaluated at $\left\langle x_{0},e_{i}\right\rangle $.

Let $\mathcal{V}_{0}$ be a fixed neighborhood of $0,$ set%
\[
\alpha_{i}=\sup_{u\in\mathcal{V}_{0}}\left\vert \frac{f_{i}\left(  u\right)
-f_{i}\left(  -u\right)  }{u\left(  f_{i}\left(  u\right)  +f_{i}\left(
-u\right)  \right)  }\right\vert
\]
and assume that :%
\[
\mathbf{A}_{4}:\sum_{i=1}^{+\infty}\alpha_{i}^{2}<+\infty.
\]

This assumption is close to those required in Mas (2007b). The next
Proposition illustrates assumption $\mathbf{A}_{4}$ in the important case when
$X$ is gaussian.

\begin{example}
Let $X$ be a centered gaussian random element in $H$ with Karhunen-Lo\`{e}ve
expansion :%
\[
X=\sum_{k=1}^{+\infty}\sqrt{\lambda_{k}}\eta_{k}e_{k}.
\]
Here the $\lambda_{k}$'s are the eigenvalues of the covariance operator of
$X,$ $\mathbb{E}\left(  X\otimes X\right)  ,$ the $e_{k}$'s are the associated
eigenvectors and the $\eta_{k}$'s are real-valued random variables $N\left(
0,1\right)  $-distributed. It is a well-known fact that $\left\langle
X,e_{k}\right\rangle =\sqrt{\lambda_{k}}\eta_{k}$ are independent real
gaussian random variables and $\mathbf{A}_{3}$ holds. Then $f_{i}\left(
u\right)  =\frac{1}{\sqrt{2\pi\lambda_{i}}}\exp\left[  -\frac{\left(
u-\left\langle x_{0},e_{i}\right\rangle \right)  ^{2}}{2\lambda_{i}}\right]  $
and%
\[
\sup_{u\in\mathcal{V}_{0}}\frac{\left\vert f_{i}\left(  u\right)
-f_{i}\left(  -u\right)  \right\vert }{u\left\vert f_{i}\left(  u\right)
+f_{i}\left(  -u\right)  \right\vert }\leq C\frac{\left\langle x_{0}%
,e_{i}\right\rangle }{\lambda_{i}}%
\]
whenever $\left\langle x_{0},e_{i}\right\rangle /\lambda_{i}\rightarrow0$ when
$i$ tends to infinity and $\mathbf{A}_{4}$ holds if :%
\begin{equation}
\sum_{i=1}^{+\infty}\frac{\left\langle x_{0},e_{i}\right\rangle ^{2}}%
{\lambda_{i}^{2}}<+\infty\label{fouq}%
\end{equation}

\end{example}

\begin{example}
We can also consider the family of densities indexed by the integer $m$ :%
\[
f_{i}\left(  u\right)  =\frac{C_{m}}{\sqrt{\lambda_{i}}}\frac{1}{1+\left(
\frac{u-\left\langle x_{0},e_{i}\right\rangle }{\sqrt{\lambda_{i}}}\right)
^{2m}}%
\]
where $C_{m}$ is a normalizing constant. We find :%
\[
\alpha_{i}\leq C\frac{\left\vert \left\langle x_{0},e_{i}\right\rangle
\right\vert }{\lambda_{i}^{m}}%
\]
and assumption $\mathbf{A}_{4}$ holds whenever the sequence $\left(
\frac{\left\vert \left\langle x_{0},e_{i}\right\rangle \right\vert }%
{\lambda_{i}^{m}}\right)  _{i\in\mathbb{N}}\in l_{2}$.
\end{example}

Since the rate of decrease of the $\lambda_{i}$'s is intimately related to the
smoothness of the random function $X,$ we may easily infer that $\mathbf{A}%
_{4}$ should be interpreted as a smoothness condition on the function $x_{0}$.
In other words, the coordinates of $x_{0}$ in the basis $e_{i}$ should tend to
zero at a rate which is significantly quicker than the eigenvalues of the
covariance operator of $X$ and hence that $x_{0}$ should be sufficiently
smoother than $X$.

It should also be noted that, when the family of densities $f_{i}$ is not
uniformly smooth enough in a neighborhood of $0,$ Assumption $\mathbf{A}_{4}$
may fail. For instance, it is not hard to see that the $\alpha_{i}$'s are not
even finite when $f_{i}$ is the density of a shifted Laplace random variable :%
\[
f_{i}\left(  u\right)  =\frac{1}{2\lambda_{i}}\exp\left(  -\frac{\left\vert
u-\left\langle x_{0},e_{i}\right\rangle \right\vert }{\lambda_{i}}\right)  .
\]

\begin{remark}
The issue of the expectation of the functional input $X$ should be raised now.
We assumed sooner that the $X_{i}$'s are centered. But in practical situations
we can expect $\mu=\mathbb{E}\left(  X\right)  $ to be a non-null function.
Then considering a new shift $x_{0}-\mu$ instead of $x_{0}$ solves the
problem. So we can always consider the centered version of $X$ but we must
take into account that any assumption made on $x_{0}$ should be valid for
$x_{0}-\mu$. For instance (\ref{fouq}) should be replaced by%
\[
\sum_{i=1}^{+\infty}\frac{\left\langle x_{0}-\mu,e_{i}\right\rangle ^{2}%
}{\lambda_{i}^{2}}<+\infty.
\]

\end{remark}

\subsection{Smoothness of the regression function}

In order to achieve our estimating procedure we cannot avoid to assume that
the function $m$ is regular. Since $m$ is a mapping from $H$ to $\mathbb{R},$
its first order derivative is an element of $\mathcal{L}\left(  H,\mathbb{R}%
\right)  ,$ the space of bounded linear functionals from $H$ to $\mathbb{R}$
which is nothing than $H^{\ast}\simeq H$. As announced sooner $m^{\prime
}\left(  x_{0}\right)  \in H.$ The second order derivative belongs to
$\mathcal{L}\left(  \mathcal{L}\left(  H,\mathbb{R}\right)  ,\mathbb{R}%
\right)  \simeq\mathcal{L}\left(  H\times H,\mathbb{R}\right)  $ and is
consequently a quadratic functional on $H\times H$ and may be represented by a
symmetric positive linear operator from $H$ to $H$ (the Hessian operator). We
will sometimes use abusive notations such as $\left\langle m^{\prime\prime
}\left(  x_{0}\right)  \left(  u\right)  ,v\right\rangle $ below and
throughout the proofs.\medskip

$\mathbf{A}_{5}:$\textit{ The first order derivative of }$m$ at $x_{0}$
$m^{\prime}\left(  x_{0}\right)  $\textit{ is defined, non null and there
exists a neighborhood }$\mathcal{V}\left(  x_{0}\right)  $ of $x_{0}$\textit{
such that :}%
\[
\mathit{\sup_{x\in\mathcal{V}\left(  x_{0}\right)  }}\left\Vert m^{\prime
\prime}\left(  x\right)  \right\Vert _{\infty}<+\infty.
\]
\medskip

This last display may be rewritten : for all $u$ in $H$ and all $x$ in a
neighborhood of $x_{0}$%
\[
\left\langle m^{\prime\prime}\left(  x\right)  u,u\right\rangle \leq
C\left\Vert u\right\Vert ^{2}.
\]

\begin{remark}
Assumption $\mathbf{A}_{5}$ assesses in a way that "the second order
derivative of $m$ in a neighborhood of $x_{0}$ is bounded".
\end{remark}

\subsection{Back to the regularized inverse}

We need for immediate purpose to define a sequence involved in the rate of
convergence of our estimate.

\begin{definition}
Let $v\left(  h\right)  $ the positive sequence defined by :%
\begin{equation}
v=v\left(  h\right)  =\left[  \mathbb{E}\left(  K\left(  \frac{\left\Vert
X-x_{0}\right\Vert }{h}\right)  \left\Vert X-x_{0}\right\Vert \rho\left(
\left\Vert X-x_{0}\right\Vert \right)  \right)  \right]  . \label{cat}%
\end{equation}

\end{definition}

It is plain that $v$ tends to zero when $h$ does.

Since they will be used in the sequel we list now some results from Mas
(2007b). They are collected in the next Proposition and consist in bounding
thre norms of operators $\Gamma_{K}$ and $\Gamma_{n,K}$

\begin{proposition}
The following bound are valid%
\begin{align}
\left\Vert \Gamma_{K}\right\Vert _{\infty}  &  \geq Cv\left(  h\right)
,\label{asyga1}\\
\left\Vert \Gamma_{n,K}-\Gamma_{K}\right\Vert _{\infty}  &  =O_{L^{2}}\left(
h^{2}\sqrt{\frac{F\left(  h\right)  }{n}}\right)  . \label{asyga2}%
\end{align}
Besides $\Gamma_{K}/v\left(  h\right)  $ may converge to a bounded operator,
say $S$, that may be compact.
\end{proposition}

Before giving the main results we have to get back to the regularized inverse
of $\Gamma_{n,K}.$ Indeed a bound on the norm of $\Gamma_{n,K}^{\dag}$ may be
derived. Under the assumption that $h^{2}F^{1/2}\left(  h\right)  /\left(
n^{1/2}v\left(  h\right)  \right)  \rightarrow0$ we see that $\left\Vert
\Gamma_{n,K}\right\Vert _{\infty}\geq Cv\left(  h\right)  $ As a consequence
of these facts we expect the norm $\Gamma_{n,K}^{\dag}$ to diverge with rate
\textit{at least} $1/v\left(  h\right)  $ since :%
\[
0<C<\left\Vert \Gamma_{n,K}\Gamma_{n,K}^{\dag}\right\Vert _{\infty}%
\leq\left\Vert \Gamma_{n,K}^{\dag}\right\Vert _{\infty}\left\Vert \Gamma
_{n,K}\right\Vert _{\infty}.
\]
If the operator $S$ mentioned in the Proposition above is compact, we may even
be aware that the norm of $v\left(  h\right)  \Gamma_{n,K}^{\dag}$ will tend
to infinity since $S^{-1}$ is unbounded whenever $S^{-1}$ exists. All this
leads us to considering the next and last assumption on $\Gamma_{n,K}^{\dag}$
:\medskip

$\mathbf{A}_{6}:$\textit{ There exists a sequence }$r_{n}\downarrow0$\textit{
such that}%
\[
\max\left\{  \left\Vert \Gamma_{K}^{\dag}\right\Vert _{\infty},\left\Vert
\Gamma_{n,K}^{\dag}\right\Vert _{\infty}\right\}  \mathit{\leq}\frac{1}%
{r_{n}v\left(  h\right)  }\quad a.s.
\]
\medskip

Here the parameter $r_{n}$ just depends on the chosen regularizing method
(penalization, Tikhonov, etc.) and may be viewed as a tuning parameter.

\begin{remark}
\label{lastrem}In fact as will be seen below the sequence $r_{n}$ may no tend
to zero. But the situation when $r_{n}\downarrow0$ is the most unfavorable one
and we intend to investigate it with care. However $r_{n}v\left(  h\right)  $
always tend to zero and cannot be bounded below because of (\ref{asyga1}) and
(\ref{asyga2}). Besides if $\Gamma_{K}/v\left(  h\right)  $ converges to an
operator with bounded inverse, the sequence $r_{n}$ can always be chosen constant.
\end{remark}

Let us take some examples to illustrate the role of $r_{n}$. We keep the
notations of display (\ref{svd}) and of the lines below.

\begin{itemize}
\item Truncated spectral regularization : remind that%
\[
\Gamma_{K}^{\dag}=\sum_{i=1}^{N_{n}}\frac{1}{\mu_{i,n}}\left(  u_{i,n}\otimes
u_{i,n}\right)
\]
where $\left(  \mu_{i,n},u_{i,n}\right)  $ are the eigenelements of
$\Gamma_{K}$ and $\left\Vert \Gamma_{K}\right\Vert _{\infty}=\sup_{i}\left\{
\mu_{i,n}\right\}  =\mu_{1,n}$ (as announced sooner the eigenvalues are
positive and arranged in a decreasing order). Hence%
\[
\left\Vert \Gamma_{K}^{\dag}\right\Vert _{\infty}=1/\inf_{1\leq i\leq N_{n}%
}\left\{  \mu_{i,n}\right\}  =\mu_{N_{n},n}%
\]
then $r_{n}=\mu_{N_{n},n}/\mu_{1,n}\downarrow0$ is the inverse of the
conditioning index of operator $\Gamma_{K}^{\dag}$.

\item Penalization : Now $\Gamma_{n,K}^{\dag}=\left(  \Gamma_{n,K}+\alpha
_{n}I\right)  ^{-1}$ with%
\[
\Gamma_{K}^{\dag}=\sum_{i=1}^{m_{n}}\frac{1}{\mu_{i,n}+\alpha_{n}}\left(
u_{i,n}\otimes u_{i,n}\right)
\]
and we can take $r_{n}=\alpha_{n}/\mu_{1,n}$. It is possible here to get
$r_{n}\uparrow+\infty$ by an accurate choice of $\alpha_{n}$ and some
information on $\mu_{1,n}.$

\item Tikhonov regularization : Here $\Gamma_{K}^{\dag}=\left(  \Gamma_{K}%
^{2}+\alpha_{n}I\right)  ^{-1}\Gamma_{K}$ and
\[
\Gamma_{K}^{\dag}=\sum_{i=1}^{m_{n}}\frac{\mu_{i,n}}{\mu_{i,n}^{2}+\alpha_{n}%
}\left(  u_{i,n}\otimes u_{i,n}\right)  .
\]
A choice for $r_{n}$ is here $\alpha_{n}/\mu_{1,n}^{2}$ and the same remark as
above holds.
\end{itemize}

\section{Statement of the results}

The central result of this article is a bound on the Mean Square Error for the
local linear estimate of the pointwise evaluation of the regression function
at a fixed design. In the sequel the generic notation $C$ stands for universal constants.

\begin{theorem}
\label{main}Fix $x_{0}$ in $H$. When assumptions $\mathbf{A}_{1}%
\mathbf{-A}_{6}$ hold and if $nF\left(  h\right)  \rightarrow+\infty$ :%
\begin{align*}
\mathbb{E}\left(  \widehat{m}_{n}\left(  x_{0}\right)  -m\left(  x_{0}\right)
\right)  ^{2}  &  \leq C\left[  \frac{h^{6}}{r_{n}^{2}}+h^{4}+\frac{h^{2}%
}{nF\left(  h\right)  }+\frac{v^{2}\left(  h\right)  }{F^{2}\left(  h\right)
}\right] \\
&  +\frac{C}{nF\left(  h\right)  }\left(  1+\frac{h^{2}}{nr_{n}v\left(
h\right)  }+\frac{v\left(  h\right)  }{r_{n}F\left(  h\right)  }\right)  .
\end{align*}
where the first line arises from the bias of our estimate and the second stems
from its variance.
\end{theorem}

\begin{remark}
If $K$ is chosen to be the naive kernel, $K\left(  s\right)  =1\!1_{\left[
0,1\right]  }\left(  s\right)  ,$ assumption $\mathbf{A}_{1}$ can be removed
and the previous theorem remains valid.
\end{remark}

\begin{remark}
It turns out that the variance term is decomposed into three. The first is
$\left(  nF\left(  h\right)  \right)  ^{-1}$ and is classical (see Ferraty,
Mas, Vieu (2007)). The two others stem directly from the underlying inverse
problem and the sequence $r_{n}$ appears.
\end{remark}

Note that we did not fix the issue of the sequence $r_{n}$ involved in the
regularizing inverses $\Gamma_{K}^{\dag}$ and $\Gamma_{n,K}^{\dag}$. Theorem
\ref{main} may be simplified under mild additional assumptions.

\begin{proposition}
\label{MA}Taking $r_{n}\asymp h$ then%
\begin{align*}
&  \mathbb{E}\left(  \widehat{m}_{n}\left(  x_{0}\right)  -m\left(
x_{0}\right)  \right)  ^{2}\\
&  \leq Ch^{4}+\frac{C}{nF\left(  h\right)  }\left(  1+\frac{h}{nv\left(
h\right)  }\right)
\end{align*}

\end{proposition}

This Proposition is derived from Theorem \ref{main} and Lemma \ref{L1}$.$

\begin{remark}
Turning back to Proposition \ref{rho} and considering displays (\ref{rho1})
and (\ref{rho2}) it is not hard to see that both functions $\rho$ are
regularly varying at $0$ with index $1+\beta$ for the first and $1$ for the
second and hence that Proposition \ref{MA} holds. It should also be noted that
from property (\ref{F2}) in Proposition \ref{bing} that we can truly expect
$\rho$ to be of index larger than $1$ whenever it is regularly varying at $0$.
This fact motivates the next Proposition.
\end{remark}

\begin{proposition}
\label{opt}Under the assumptions of Theorem \ref{main} and of Proposition
\ref{MA}, if the auxiliary function $\rho$ is regularly varying at $0$ with
index $g\geq1$,%
\begin{equation}
v\left(  h\right)  \asymp h\rho\left(  h\right)  F\left(  h\right)  .
\label{MC}%
\end{equation}
Then if $\rho\left(  s\right)  \geq Cs^{4}$ in a neighborhood of $0,$ the mean
square error becomes :%
\[
\mathbb{E}\left(  \widehat{m}_{n}\left(  x_{0}\right)  -m\left(  x_{0}\right)
\right)  ^{2}\leq C\left(  h^{4}+\frac{1}{nF\left(  h\right)  }\right)
\]
and the rate of decrease of the Mean Square Error depends on $h^{\ast}$ given
by%
\begin{equation}
\left(  h^{\ast}\right)  ^{4}F\left(  h^{\ast}\right)  =\frac{1}{n}.
\label{makeba}%
\end{equation}
If $\rho\left(  s\right)  /s^{4}\rightarrow0$ when $s\rightarrow0$ the above
rate is damaged$.$ For instance taking $r_{n}\asymp h$ the MSE becomes :%
\[
\mathbb{E}\left(  \widehat{m}_{n}\left(  x_{0}\right)  -m\left(  x_{0}\right)
\right)  ^{2}\leq C\left(  h^{4}+\frac{1}{n^{2}F^{2}\left(  h\right)
\rho\left(  h\right)  }\right)  .
\]

\end{proposition}

\begin{remark}
Display \ref{MC} was proved in Mas (2007b). In the first case (when
$\rho\left(  s\right)  \geq Cs^{4}$), since the bias term is here an $O\left(
h^{4}\right)  $, the rate of convergence of our estimate outperfoms the one
computed in Ferraty, Mas, Vieu (2007). The estimate was a classical
Nadaraya-Watson kernel estimator whose bias was an $O\left(  h^{2}\right)  .$
Obviously the rate of convergence in the second case is damaged but even for
very irregular processes such as Brownian motion or Brownian Bridge function
$\rho\left(  s\right)  $ is above $s^{2}$ or $s^{3}$ depending on the norms
that are used. The interested reader is referred for instance to displays (20)
and (22) in Mayer-Wolf, Zeitouni (1993) or Proposition 6.1 p.568 in Li, Shao
(2002) but will have to carry out some additional computations. It seems
reasonable to think that this unfavorable situation will rarely occur in a
usual statistical context (with functions reconstructed on "smooth spaces").
However we prove just below that, even when $\rho$ decays rapidly to $0,$ it
is always possible to choose a regularizing method for $\Gamma_{n,K}$ that
reaches the best rate of display (\ref{makeba}).
\end{remark}

\begin{remark}
It may be fruitful for practical purposes to comment on formula (\ref{makeba}%
). First we see that when $X\in\mathbb{R}^{d},$ $F\left(  h\right)  \sim
Ch^{d}$ then the rate of convergence in mean square turns out to be
$n^{-2/\left(  4+d\right)  }$ which is the optimal rate of convergence for a
twice-differentiable regression function (see Stone (1982)). When the small
ball probability belongs to the class $\Gamma_{0},$ this rate will depend on
$\rho$. We know that the term $F\left(  h\right)  $ will always tend to $0$
quicker than $h^{4}$ and will consequently determine the choice of $h$. The
situation is consequently more intricate than in the multivariate setting.
However following the example of displays (\ref{sbp1}) and (\ref{sbp2}) we get
repectively%
\begin{align*}
h_{n}^{\ast}  &  =C\left(  \log n\right)  ^{-1/\beta}\\
h_{n}^{\ast}  &  =C\left(  \log n\right)  ^{-1/2}%
\end{align*}
where $\beta\leq3$ when $\rho\left(  s\right)  \geq Cs^{4}$. Finally the rate
of decrease of the mean square error is a $O\left(  \left(  \log n\right)
^{-c}\right)  $ where $c>1.$
\end{remark}

The last Proposition is devoted to dealing with the situation described along
Remark \ref{lastrem} : when $r_{n}$ does not tend to zero. This cannot happen
when the regularizing method is the spectral truncation but may occur when
either a penalization or a Tikhonov method are applied. We remind that we
cannot avoid the condition $r_{n}v\left(  h\right)  \downarrow0$. We start
from Theorem \ref{main}.

\begin{proposition}
When assumptions $\mathbf{A}_{1}\mathbf{-A}_{6}$ hold, if $nF\left(  h\right)
\rightarrow+\infty$, when the regularizing method allows to do so, taking
$r\left(  h\right)  =1/\rho\left(  h\right)  $ provides :%
\[
\mathbb{E}\left(  \widehat{m}_{n}\left(  x_{0}\right)  -m\left(  x_{0}\right)
\right)  ^{2}\leq Ch^{4}+C\frac{1}{nF\left(  h\right)  }.
\]

\end{proposition}

Obviously $r_{n}v\left(  h\right)  $ tends to $0$. If the chosen method is
penalization such that $\Gamma_{n,K}^{\dag}=\left(  \Gamma_{n,K}+\alpha
_{n}S\right)  ^{-1}$ it suffices to take $\alpha_{n}=h^{\ast}F\left(  h^{\ast
}\right)  $ to achieve our goal. The proof of this Proposition is easy hence omitted.

\begin{remark}
The rate obtained at display (\ref{makeba}) issued from Proposition \ref{opt}
should be compared with the minimax rate obtained by Fan (1993) for scalar
inputs. The MSE was then $Ch^{4}+C/\left(  nh\right)  .$ We see that,
replacing $F\left(  h\right)  $ by $h$ (which is logic if we consider the
remark about the multivariate case just below display (\ref{sbp1}) in the
section devoted to the small ball problems), both formulas match. This fact
leads us to another interesting issue : does this rate inherit the optimal
(minimax) properties found by Fan in his article ? This question goes beyond
the scope of this article. Besides not much has been done until now about
optimal estimation for functional data -to the authors' knowledge. But there
is no doubt that this issue will be addressed in the next future.
\end{remark}

\section{Conclusion}

Obviously this article could be the starting point for other issues such as
almost sure or weak convergence of the estimate. Almost all practical aspects
were left out on purpose : they will certainly give birth to another article.
However the main goal of this essentially theoretic work was to underline the
rather large scope of our study. We had to seek several ideas in such various
areas as probability theory, functional analysis, statistical theory of
extremes, inverse problems theory. Finally it turns out that it is possible to
get, in the functional setting, almost the same rate of decay for the bias as
in the case of scalar inputs. The variance involves the small ball probability
evaluated at $h,$ the selected bandwidth. A drawback arises with the necessity
to introduce a new parameter : the regularizing sequence $r_{n},$ which
depends on the sample size (more precisely on the bandwidth $h$). We give no
clue to find out in practical situations the bandwidth $h$ but we guess that
the ever wider literature on functional data will quickly overcome this
problem by adapting classical methods such as cross-validation for instance.

Another major practical concern relies in the estimation of the unknown
auxiliary function $\rho$. Several tracks already appear to address this
issue. One may think of adapting some techniques from extreme theory. After
all $\rho$ characterizes the extreme behaviour of $\left\Vert X\right\Vert $
like tail indices for Weibull or Pareto distributions. The only difference
stems from the fact that $\rho$ is a function and not just a real number. The
other idea lies in the article by Mas (2008) where the auxiliary function
$\rho$ is explicitely linked with the eigenvalues of the ordinary covariance
operator of $X.$ From the estimation of these eigenvalues (which is a basic
procedure) it should be possible to propose a consistent estimation of the
auxiliary function as a by-product.

\section{Proofs}

For the sake of clarity we begin with an outline of the proofs. The following
bias-variance decomposition for $\widehat{m}_{n}\left(  x_{0}\right)
-m\left(  x_{0}\right)  $ holds :%
\begin{align*}
\widehat{m}_{n}\left(  x_{0}\right)  -m\left(  x_{0}\right)   &  =\dfrac
{\sum_{i=1}^{n}y_{i}\omega_{i,n}}{\sum_{i=1}^{n}\omega_{i,n}}-m\left(
x_{0}\right) \\
&  =\dfrac{\sum_{i=1}^{n}\left(  y_{i}-m\left(  x_{0}\right)  \right)
\omega_{i,n}}{\sum_{i=1}^{n}\omega_{i,n}}\\
&  =\dfrac{\sum_{i=1}^{n}\left(  y_{i}-m\left(  X_{i}\right)  \right)
\omega_{i,n}}{\sum_{i=1}^{n}\omega_{i,n}}+\dfrac{\sum_{i=1}^{n}\left(
m\left(  X_{i}\right)  -m\left(  x_{0}\right)  \right)  \omega_{i,n}}%
{\sum_{i=1}^{n}\omega_{i,n}}.
\end{align*}
We denote :%
\begin{align}
T_{b,n}  &  =\dfrac{\sum_{i=1}^{n}\left(  m\left(  X_{i}\right)  -m\left(
x_{0}\right)  \right)  \omega_{i,n}}{\sum_{i=1}^{n}\omega_{i,n}}%
,\label{bias}\\
T_{v,n}  &  =\dfrac{\sum_{i=1}^{n}\left(  y_{i}-m\left(  X_{i}\right)
\right)  \omega_{i,n}}{\sum_{i=1}^{n}\omega_{i,n}}\label{variance}\\
&  =\dfrac{\sum_{i=1}^{n}\omega_{i,n}\varepsilon_{i}}{\sum_{i=1}^{n}%
\omega_{i,n}}\nonumber
\end{align}
where $\varepsilon$ was defined at display (\ref{modele}). Here $T_{b,n}$ is a
bias term and $T_{v,n}$ is a variance term. Finally we get :%
\begin{equation}
\mathbb{E}\left[  \widehat{m}_{n}\left(  x_{0}\right)  -m\left(  x_{0}\right)
\right]  ^{2}=\mathbb{E}T_{b,n}^{2}+\mathbb{E}T_{v,n}^{2}+2\mathbb{E}\left(
T_{b,n}T_{v,n}\right)  \label{moon}%
\end{equation}
and since%
\begin{align*}
\mathbb{E}\left(  T_{b,n}T_{v,n}\right)   &  =\mathbb{E}\left(  T_{b,n}%
\mathbb{E}\left(  T_{v,n}|X_{1},...,X_{n}\right)  \right) \\
&  =0
\end{align*}
computing the mean square error of $\widehat{m}_{n}\left(  x_{0}\right)  $
comes down to computing $\mathbb{E}T_{b,n}^{2}$ and $\mathbb{E}T_{v,n}^{2}$
which will be done later.

The proof section is tiled into two subsections. The first one is devoted to
giving preliminary results as well as Lemmas. In the second the main results
are derived.

\subsection{Preliminary results}

We assume that assumptions $\mathbf{A}_{1}\mathbf{-A}_{6}$ hold once and for
all. The next two Lemmas are given for further purposes. Their proofs are
omitted. The interested reader will find them in Mas (2007b).

\begin{lemma}
\label{LP1}If $f$ belongs to the class $\Gamma_{0}$ with auxiliary function
$\rho$, then for all $p\in\mathbb{N}$,%
\[
\int_{0}^{1}\frac{t^{p}}{\sqrt{1-t^{2}}}f\left(  s\sqrt{1-t^{2}}\right)
dt\underset{s\rightarrow0}{\sim}2^{\frac{p-1}{2}}\Gamma\left(  \frac{p+1}%
{2}\right)  f\left(  s\right)  \left(  \frac{\rho\left(  s\right)  }%
{s}\right)  ^{\frac{p+1}{2}}.
\]

\end{lemma}

For any $x=\sum x_{k}e_{k}$ in $H$ and for $i\in\mathbb{N}$ set $\left\Vert
x\right\Vert _{\neq i}^{2}=\sum_{k\neq i}x_{k}^{2}.$

We denote $f_{\neq i}$ the density of $\left\Vert X-x_{0}\right\Vert _{\neq
i}$. We need to compute both densities $f_{\left\Vert X-x_{0}\right\Vert }$
(density of $\left\Vert X-x_{0}\right\Vert $) and $f_{\left\langle
X-x_{0},e_{i}\right\rangle ,\left\Vert X-x_{0}\right\Vert }$ (density of the
couple $\left(  \left\langle X-x_{0},e_{i}\right\rangle ,\left\Vert
X-x_{0}\right\Vert \right)  $).

\begin{lemma}
\label{LP2}We have :%
\begin{align}
f_{\left\langle X-x_{0},e_{i}\right\rangle ,\left\Vert X-x_{0}\right\Vert
}\left(  u,v\right)   &  =\frac{v}{\sqrt{v^{2}-u^{2}}}f_{i}\left(  u\right)
f_{\neq i}\left(  \sqrt{v^{2}-u^{2}}\right)  1\hspace{-3pt}1_{\left\{
v\geq\left\vert u\right\vert \right\}  },\label{d1}\\
f_{\left\Vert X-x_{0}\right\Vert }\left(  v\right)   &  =v\int_{-1}^{1}%
\frac{f_{i}\left(  vt\right)  }{\sqrt{1-t^{2}}}f_{\neq i}\left(
v\sqrt{1-t^{2}}\right)  dt. \label{d2}%
\end{align}
Besides if $f_{\left\Vert X-x_{0}\right\Vert }$ and $f_{\neq i}$ are $\Gamma
$-varying for all $i$ then they have all $\rho$ as auxiliary function.
\end{lemma}

We begin with more specific computational\ Lemmas.

\begin{lemma}
\label{L1}Let $\varphi$ be a positive real valued function, bounded on
$\left[  0,1\right]  $ and regularly varying at $0$ with index $g\geq1$ and
let $p\in\mathbb{N}$ :%
\begin{equation}
\mathbb{E}K^{p}\left(  \left\Vert \frac{X-x_{0}}{h}\right\Vert \right)
\varphi\left(  \left\Vert X-x_{0}\right\Vert \right)  \underset{h\rightarrow
0}{\sim}K^{p}\left(  1\right)  \varphi\left(  h\right)  F\left(  h\right)  .
\label{sting}%
\end{equation}
As important special cases we mention :%
\begin{align*}
\mathbb{E}K\left(  \left\Vert \frac{X-x_{0}}{h}\right\Vert \right)   &  \sim
K\left(  1\right)  F\left(  h\right)  ,\quad\mathbb{E}K^{2}\left(  \left\Vert
\frac{X-x_{0}}{h}\right\Vert \right)  \sim K^{2}\left(  1\right)  F\left(
h\right)  ,\\
\mathbb{E}\left[  \left\Vert X-x_{0}\right\Vert ^{m}K\left(  \left\Vert
\frac{X-x_{0}}{h}\right\Vert \right)  \right]   &  \sim K\left(  1\right)
F\left(  h\right)  h^{m}.
\end{align*}

\end{lemma}

\textbf{Proof :}

We prove (\ref{sting}) when $p=1$ and denote $\mathbb{P}^{\left\Vert
X_{i}-x_{0}\right\Vert /h}$ the distribution of the random variable
$\left\Vert X_{i}-x_{0}\right\Vert /h$. Since%
\[
\mathbb{E}K\left(  \left\Vert \frac{X-x_{0}}{h}\right\Vert \right)
\varphi\left(  \left\Vert X-x_{0}\right\Vert \right)  =\int_{0}^{1}K\left(
u\right)  \varphi\left(  hu\right)  d\mathbb{P}^{\left\Vert X_{i}%
-x_{0}\right\Vert /h}\left(  u\right)  ,
\]
and from $K\left(  u\right)  \varphi\left(  hu\right)  =K\left(  1\right)
\varphi\left(  h\right)  -\int_{u}^{1}\left[  K\left(  s\right)
\varphi\left(  hs\right)  \right]  ^{\prime}ds$ we get :%
\begin{align*}
&  \mathbb{E}K\left(  \left\Vert \frac{X-x_{0}}{h}\right\Vert \right)
\varphi\left(  \left\Vert X-x_{0}\right\Vert \right) \\
&  =K\left(  1\right)  \varphi\left(  h\right)  \int d\mathbb{P}^{\left\Vert
X_{i}-x_{0}\right\Vert /h}\left(  u\right)  -\int\int_{0\leq u\leq s\leq
1}\left[  K\left(  s\right)  \varphi\left(  hs\right)  \right]  ^{\prime
}d\mathbb{P}^{\left\Vert X_{i}-x_{0}\right\Vert /h}\left(  u\right)
\end{align*}
Applying Fubini's Theorem we get :%
\begin{align*}
\mathbb{E}K\left(  \left\Vert \frac{X-x_{0}}{h}\right\Vert \right)
\varphi\left(  \left\Vert X-x_{0}\right\Vert \right)   &  =K\left(  1\right)
\varphi\left(  h\right)  F\left(  h\right)  -\int\left[  K\left(  s\right)
\varphi\left(  hs\right)  \right]  ^{\prime}F\left(  hs\right)  ds\\
&  =K\left(  1\right)  \varphi\left(  h\right)  F\left(  h\right)  \left(
1-\mathcal{R}_{h}\right)
\end{align*}
with%
\[
\mathcal{R}_{h}=\int_{0}^{1}\frac{K^{\prime}\left(  s\right)  \varphi\left(
hs\right)  +K\left(  s\right)  h\varphi^{\prime}\left(  hs\right)  }%
{\varphi\left(  h\right)  }\frac{F\left(  hs\right)  }{F\left(  h\right)  }ds
\]
Since $F$ is gamma-varying at $0,$ display (\ref{F1b}) in Proposition
\ref{bing} tells us that $F\left(  hs\right)  /F\left(  h\right)
\rightarrow0$ as $h\rightarrow0$. As $\varphi$ is regularly varying at $0$
with, say, index $g\geq1,$ $\varphi\left(  hs\right)  /\varphi\left(
h\right)  \rightarrow s^{g}$ as $h$ goes to $0$. Remind also that $K^{\prime}$
is integrable. We deal with
\[
h\frac{\varphi^{\prime}\left(  hs\right)  }{\varphi\left(  h\right)  }%
=h\frac{\varphi^{\prime}\left(  h\right)  }{\varphi\left(  h\right)  }%
\frac{\varphi^{\prime}\left(  hs\right)  }{\varphi^{\prime}\left(  h\right)  }%
\]
Now in Bingham, Goldie and Teugels (1987), the definition of regular variation
is given p.18. From Theorem 1.7.2b p.39 we deduce that $\varphi^{\prime}$ is
regularly varying with index $g-1\geq0$ hence that :%
\[
\lim_{h\rightarrow0}\frac{\varphi^{\prime}\left(  hs\right)  }{\varphi
^{\prime}\left(  h\right)  }=s^{g-1}%
\]
uniformly with respect to $s\in\left]  0,1\right]  $ and by the direct part of
Karamata's Theorem p.28 (take $g=\rho-1,$ $\sigma=0$ and $f=\varphi^{\prime}$)
that :%
\[
\limsup_{h\rightarrow0}h\frac{\varphi^{\prime}\left(  h\right)  }%
{\varphi\left(  h\right)  }\rightarrow g
\]
which means that $h\varphi^{\prime}\left(  hs\right)  /\varphi\left(
h\right)  $ converges pointwise to $gs^{g-1}$ (which is integrable with
respect to Lebesgue's measure). Then we can apply Lebesgue's dominated
convergence theorem and Proposition \ref{bing} (see display (\ref{F1b})) to
get $\mathcal{R}_{h}\rightarrow0$ as $h\rightarrow0$. This last step leads to
the announced result.\medskip

For the sake of shortness we will sometimes set :%
\[
Z=X-x_{0},K=K\left(  \left\Vert X-x_{0}\right\Vert /h\right)
\]
and :%
\[
\overline{Z}_{K,n}=\frac{1}{n}\sum_{i=1}^{n}Z_{i}K_{i}=\frac{1}{n}\sum
_{k=1}^{n}\left(  X_{i}-x_{0}\right)  K\left(  \left\Vert X_{i}-x_{0}%
\right\Vert /h\right)  .
\]

The next lemma is a crucial.

\begin{lemma}
\label{L2}We have :%
\begin{align*}
\left\Vert \mathbb{E}\left[  ZK\right]  \right\Vert ^{2}  &  =\left\Vert
\mathbb{E}\left[  K\left(  \left\Vert \frac{X-x_{0}}{h}\right\Vert \right)
\left(  X-x_{0}\right)  \right]  \right\Vert ^{2}\\
&  \leq Cv^{2}\left(  h\right)  .
\end{align*}

\end{lemma}

\begin{remark}
We can measure the sharpness of the previous bound. Indeed a very simple
inequality would give by\ Lemma \ref{L1} :%
\[
\left\Vert \mathbb{E}\left[  ZK\right]  \right\Vert ^{2}\leq\left(
\mathbb{E}\left\Vert ZK\right\Vert \right)  ^{2}=\left(  \mathbb{E}\left\Vert
Z\right\Vert K\right)  ^{2}\sim Ch^{2}F^{2}\left(  h\right)
\]
whereas in view of (\ref{cat}) and -when $\rho$ is regularly varying at $0$
with positive index- of Lemma \ref{L1},%
\[
\left[  \mathbb{E}\left(  K\left(  \frac{\left\Vert X-x_{0}\right\Vert }%
{h}\right)  \left\Vert X-x_{0}\right\Vert \rho\left(  \left\Vert
X-x_{0}\right\Vert \right)  \right)  \right]  ^{2}\leq h^{2}\rho^{2}\left(
h\right)  F^{2}\left(  h\right)  .
\]
So the bound was improved by a rate of $\rho^{2}\left(  h\right)  =o\left(
h^{2}\right)  .$
\end{remark}

\textbf{Proof :}

Computations here are quite similar but however distinct from those carried
out in Mas (2007b). We start with projecting $\mathbb{E}\left[  K\left(
\left\Vert \frac{X-x_{0}}{h}\right\Vert \right)  \left(  X-x_{0}\right)
\right]  $ on the basis $\left(  e_{i}\right)  _{i\in\mathbb{N}}$ mentioned in
$\mathbf{A}_{3}$ and compute thanks to Lemma \ref{LP2} :%
\begin{align}
&  \mathbb{E}\left[  K\left(  \frac{\left\Vert X-x_{0}\right\Vert }{h}\right)
\left\langle X-x_{0},e_{i}\right\rangle \right] \nonumber\\
&  =\int\int K\left(  \frac{v}{h}\right)  \frac{uv}{\sqrt{v^{2}-u^{2}}}%
f_{i}\left(  u\right)  f_{\neq i}\left(  \sqrt{v^{2}-u^{2}}\right)
1\hspace{-3pt}1_{\left\{  h\geq v\geq\left\vert u\right\vert \right\}
}dudv\nonumber\\
&  =\int_{0}^{h}vK\left(  \frac{v}{h}\right)  \left(  \int_{-v}^{v}\frac
{u}{\sqrt{v^{2}-u^{2}}}f_{i}\left(  u\right)  f_{\neq i}\left(  \sqrt
{v^{2}-u^{2}}\right)  du\right)  dv.\label{Malcolm}\\
& \nonumber
\end{align}

Now we deal with%
\begin{align*}
&  \int_{-v}^{v}\frac{u}{\sqrt{v^{2}-u^{2}}}f_{i}\left(  u\right)  f_{\neq
i}\left(  \sqrt{v^{2}-u^{2}}\right)  du\\
&  =\left(  \int_{0}^{v}\frac{u}{\sqrt{v^{2}-u^{2}}}\left(  f_{i}\left(
u\right)  -f_{i}\left(  -u\right)  \right)  f_{\neq i}\left(  \sqrt
{v^{2}-u^{2}}\right)  du\right)
\end{align*}
hence%
\begin{align*}
&  \left\vert \int_{-v}^{v}\frac{u}{\sqrt{v^{2}-u^{2}}}f_{i}\left(  u\right)
f_{\neq i}\left(  \sqrt{v^{2}-u^{2}}\right)  du\right\vert \\
&  \leq\sup_{0\leq u\leq v\leq h}\left\vert \frac{1}{u}\frac{f_{i}\left(
u\right)  -f_{i}\left(  -u\right)  }{f_{i}\left(  u\right)  +f_{i}\left(
-u\right)  }\right\vert \int_{0}^{v}\frac{u^{2}}{\sqrt{v^{2}-u^{2}}}\left(
f_{i}\left(  u\right)  +f_{i}\left(  -u\right)  \right)  f_{\neq i}\left(
\sqrt{v^{2}-u^{2}}\right)  du\\
&  \leq\alpha_{i}\int_{-v}^{v}\frac{u^{2}}{\sqrt{v^{2}-u^{2}}}f_{i}\left(
u\right)  f_{\neq i}\left(  \sqrt{v^{2}-u^{2}}\right)  du.
\end{align*}
As a consequence of the preceding lines we get%
\[
\left\vert \mathbb{E}\left[  K\left(  \frac{\left\Vert X-x_{0}\right\Vert }%
{h}\right)  \left\langle X-x_{0},e_{i}\right\rangle \right]  \right\vert
\leq\alpha_{i}\mathbb{E}\left[  K\left(  \frac{\left\Vert X-x_{0}\right\Vert
}{h}\right)  \left\langle X-x_{0},e_{i}\right\rangle ^{2}\right]
\]
leading to%
\begin{align*}
\left\Vert \mathbb{E}\left[  K\left(  \frac{\left\Vert X-x_{0}\right\Vert }%
{h}\right)  \left(  X-x_{0}\right)  \right]  \right\Vert ^{2}  &  \leq
\sum_{i=1}^{+\infty}\alpha_{i}^{2}\left(  \mathbb{E}\left[  K\left(
\frac{\left\Vert X-x_{0}\right\Vert }{h}\right)  \left\langle X-x_{0}%
,e_{i}\right\rangle ^{2}\right]  \right)  ^{2}\\
&  \leq\sup_{i}\left(  \mathbb{E}\left[  K\left(  \frac{\left\Vert
X-x_{0}\right\Vert }{h}\right)  \left\langle X-x_{0},e_{i}\right\rangle
^{2}\right]  \right)  ^{2}\sum_{i=1}^{+\infty}\alpha_{i}^{2}\\
&  \leq Cv\left(  h\right)  ^{2}.
\end{align*}

\begin{lemma}
\label{L3}Both following bounds hold :%
\begin{align*}
\mathbb{E}\left\Vert \overline{Z}_{K,n}-\mathbb{E}\overline{Z}_{K,n}%
\right\Vert ^{2}  &  \leq C\frac{h^{2}F\left(  h\right)  }{n},\\
\mathbb{E}\left\Vert \overline{Z}_{K,n}-\mathbb{E}\overline{Z}_{K,n}%
\right\Vert ^{4}  &  \leq C\frac{h^{4}F^{2}\left(  h\right)  }{n^{2}}.
\end{align*}

\end{lemma}

\textbf{Proof : }We may see $\overline{Z}_{K,n}-\mathbb{E}\overline{Z}_{K,n}$
as an array of $n$ independent centered random element with values in a
Hilbert space. Denote :%
\[
U=ZK-\mathbb{E}\left(  ZK\right)  .
\]
Then $\overline{Z}_{K,n}-\mathbb{E}\overline{Z}_{K,n}=\left(  1/n\right)
\sum_{k=1}^{n}U_{k}$. We limit ourselves to proving the second display in the
Lemma, which is the most technical. It is a slightly painful but however quite
simple calculation to get :%
\begin{align*}
\mathbb{E}\left\Vert \overline{Z}_{K,n}-\mathbb{E}\overline{Z}_{K,n}%
\right\Vert ^{4}  &  \leq C\left[  \frac{1}{n^{3}}\mathbb{E}\left\Vert
U_{1}\right\Vert ^{4}+\frac{1}{n^{2}}\left(  \mathbb{E}\left\Vert
U_{1}\right\Vert ^{2}\right)  ^{2}+\frac{1}{n^{2}}\mathbb{E}\left\langle
U_{1},U_{2}\right\rangle \left\langle U_{1},U_{2}\right\rangle \right] \\
&  \leq C\left[  \frac{1}{n^{3}}\mathbb{E}\left\Vert U_{1}\right\Vert
^{4}+\frac{2}{n^{2}}\left(  \mathbb{E}\left\Vert U_{1}\right\Vert ^{2}\right)
^{2}\right]
\end{align*}
where the last line stems from the first by Cauchy-Schwarz inequality. We do
not want to go too deeply into steps that may be easily deduced and we hope
the reader will agree that, due to the $n^{3}$ denominator the first term on
the right in the display above may be neglected with respect to the second
(even if $\left(  \mathbb{E}\left\Vert U_{1}\right\Vert ^{2}\right)  ^{2}%
\leq\mathbb{E}\left\Vert U_{1}\right\Vert ^{4}$). We turn to :%
\[
\mathbb{E}\left\Vert U\right\Vert ^{2}=\mathbb{E}\left[  \left\Vert
Z\right\Vert ^{2}K^{2}\right]  -\left\Vert \mathbb{E}\left[  ZK\right]
\right\Vert ^{2}.
\]
It follows from Lemma \ref{L1} and Lemma \ref{L2} that :%
\[
\left\Vert \mathbb{E}\left[  ZK\right]  \right\Vert ^{2}=o\left(
\mathbb{E}\left[  \left\Vert Z\right\Vert ^{2}K^{2}\right]  \right)
\]
hence that :%
\begin{align*}
\mathbb{E}\left\Vert U\right\Vert ^{2}  &  \sim\mathbb{E}\left[  \left\Vert
Z\right\Vert ^{2}K^{2}\right] \\
&  \sim Ch^{2}F\left(  h\right)
\end{align*}
which finishes the proof of Lemma \ref{L3}.

\begin{lemma}
\label{L4}We have :%
\[
\mathbb{E}\omega_{1,n}^{2}\leq C\left(  F\left(  h\right)  +\frac
{h^{2}F\left(  h\right)  }{nr_{n}v\left(  h\right)  }+\frac{v\left(  h\right)
}{r_{n}}\right)  .
\]

\end{lemma}

\textbf{Proof : }Developping $\omega_{1,n}^{2}$ we get :%
\begin{align*}
\omega_{1,n}^{2}  &  =K_{1}^{2}+\left\langle \Gamma_{n,K}^{\dag}\overline
{Z}_{K,n},K_{1}Z_{_{1}}\right\rangle ^{2}-2K_{1}^{2}\left\langle \Gamma
_{n,K}^{\dag}\overline{Z}_{K,n},Z_{_{1}}\right\rangle \\
&  \leq2K_{1}^{2}+2\left\langle \Gamma_{n,K}^{\dag}\overline{Z}_{K,n}%
,K_{1}Z_{_{1}}\right\rangle ^{2}.
\end{align*}
We deal essentially with the second term since by\ Lemma \ref{L1} we know that
$\mathbb{E}K_{1}^{2}=O\left(  F\left(  h\right)  \right)  .$ We have~:%
\[
\mathbb{E}\left\langle \Gamma_{n,K}^{\dag}\overline{Z}_{K,n},K_{1}Z_{_{1}%
}\right\rangle ^{2}\leq C\mathbb{E}\left\langle \Gamma_{n,K}^{\dag}%
\overline{Z}_{K,n},\sqrt{K_{1}}Z_{_{1}}\right\rangle ^{2}%
\]
where $C$ is here nothing but $\sup_{s}\left\vert \sqrt{K\left(  s\right)
}\right\vert $. Since the expectation in the above display we bay rewritten :%
\[
\mathbb{E}\left\langle \Gamma_{n,K}^{\dag}\overline{Z}_{K,n},\sqrt{K_{i}}%
Z_{i}\right\rangle ^{2}%
\]
for all $i$ we also have :%
\begin{align*}
\mathbb{E}\left\langle \Gamma_{n,K}^{\dag}\overline{Z}_{K,n},K_{1}Z_{_{1}%
}\right\rangle ^{2}  &  \leq\frac{C}{n}\sum_{i=1}^{n}\mathbb{E}\left\langle
\Gamma_{n,K}^{\dag}\overline{Z}_{K,n},\sqrt{K_{i}}Z_{i}\right\rangle ^{2}\\
&  =C\mathbb{E}\left[  \frac{1}{n}\sum_{i=1}^{n}K_{i}\left\langle \Gamma
_{n,K}^{\dag}\overline{Z}_{K,n},Z_{i}\right\rangle \left\langle \Gamma
_{n,K}^{\dag}\overline{Z}_{K,n},Z_{i}\right\rangle \right] \\
&  =C\mathbb{E}\left(  \left\langle \Gamma_{n,K}\left(  \Gamma_{n,K}^{\dag
}\overline{Z}_{K,n}\right)  ,\Gamma_{n,K}^{\dag}\overline{Z}_{K,n}%
\right\rangle \right)
\end{align*}
since for all $u$ in $H$%
\[
\frac{1}{n}\sum_{i=1}^{n}K_{i}\left\langle u,Z_{i}\right\rangle \left\langle
u,Z_{i}\right\rangle =\left\langle \Gamma_{n,K}\left(  u\right)
,u\right\rangle .
\]
At last%
\[
\mathbb{E}\left\langle \Gamma_{n,K}^{\dag}\overline{Z}_{K,n},K_{1}Z_{_{1}%
}\right\rangle ^{2}\leq C\mathbb{E}\left(  \left\langle \Gamma_{n,K}^{\dag
}\Gamma_{n,K}\Gamma_{n,K}^{\dag}\overline{Z}_{K,n},\overline{Z}_{K,n}%
\right\rangle \right)  .
\]
We set $S_{n}=\Gamma_{n,K}^{\dag}\Gamma_{n,K}\Gamma_{n,K}^{\dag}$ ($S_{n}$ is
a positive symmetic operator) and notice that :%
\[
\left\Vert S_{n}\right\Vert _{\infty}\leq C\left\Vert \Gamma_{n,K}^{\dag
}\right\Vert _{\infty}\leq\frac{C}{r_{n}v\left(  h\right)  }\quad a.s.
\]
because $\sup_{n}\left\Vert \Gamma_{n,K}\Gamma_{n,K}^{\dag}\right\Vert
_{\infty}<+\infty.$

Our last inequality becomes :%
\begin{align*}
\mathbb{E}\left\langle \Gamma_{n,K}^{\dag}\overline{Z}_{K,n},K_{1}Z_{_{1}%
}\right\rangle ^{2}  &  \leq C\mathbb{E}\left(  \left\langle S_{n}\overline
{Z}_{K,n},\overline{Z}_{K,n}\right\rangle \right) \\
&  =C\mathbb{E}\left(  \left\Vert S_{n}^{1/2}\overline{Z}_{K,n}\right\Vert
^{2}\right) \\
&  \leq C\left(  \mathbb{E}\left\Vert S_{n}^{1/2}\left(  \overline{Z}%
_{K,n}-\mathbb{E}\overline{Z}_{K,n}\right)  \right\Vert ^{2}+\mathbb{E}%
\left\Vert S_{n}^{1/2}\left(  \mathbb{E}\overline{Z}_{K,n}\right)  \right\Vert
^{2}\right) \\
&  \leq\frac{C}{r_{n}v\left(  h\right)  }\mathbb{E}\left\Vert \left(
\overline{Z}_{K,n}-\mathbb{E}\overline{Z}_{K,n}\right)  \right\Vert ^{2}%
+\frac{C}{r_{n}v\left(  h\right)  }\left\Vert \mathbb{E}\left(  ZK\right)
\right\Vert ^{2}.
\end{align*}

We invoke Lemma \ref{L3} and Lemma \ref{L2} to bound both terms in the
preceding display. At last we get :%
\[
\mathbb{E}\left\langle \Gamma_{n,K}^{\dag}\overline{Z}_{K,n},K_{1}Z_{_{1}%
}\right\rangle ^{2}\leq C\left(  \frac{h^{2}F\left(  h\right)  }%
{nr_{n}v\left(  h\right)  }+\frac{v\left(  h\right)  }{r_{n}}\right)
\]
which yields the desired result.

\begin{lemma}
\label{L5}When $nF\left(  h\right)  \rightarrow+\infty,$%
\[
\frac{\sum_{i=1}^{n}K_{i}}{nK\left(  1\right)  F\left(  h\right)  }%
-1\overset{L^{2}}{\rightarrow}0.
\]

\end{lemma}

where $\overset{L^{2}}{\rightarrow}$ denotes convergence in mean square.

\textbf{Proof : }%
\[
\frac{\sum_{i=1}^{n}K_{i}}{nK\left(  1\right)  F\left(  h\right)  }%
-1=\frac{\sum_{i=1}^{n}\left(  K_{i}-\mathbb{E}K_{i}\right)  }{nK\left(
1\right)  F\left(  h\right)  }+\frac{\mathbb{E}K_{i}}{K\left(  1\right)
F\left(  h\right)  }-1.
\]
By Lemma \ref{L1} the second term tends to zero. We deal with the first one.
We note that :%
\[
\mathbb{E}\left(  K_{i}-\mathbb{E}K_{i}\right)  ^{2}=\mathbb{E}K_{i}%
^{2}-\left(  \mathbb{E}K_{i}\right)  ^{2}\sim K^{2}\left(  1\right)  F\left(
h\right)
\]
by Lemma \ref{L1} again. Straightforward computations give :%
\[
\frac{\sum_{i=1}^{n}\left(  K_{i}-\mathbb{E}K_{i}\right)  }{nK\left(
1\right)  F\left(  h\right)  }=O_{L^{2}}\left(  \frac{1}{\sqrt{nF\left(
h\right)  }}\right)
\]
hence the conclusion.

\begin{lemma}
\label{L6}We have :%
\[
\mathbb{E}\left\langle \Gamma_{n,K}^{\dag}\overline{Z}_{K,n},\overline
{Z}_{K,n}\right\rangle ^{2}\leq C\frac{h^{4}F^{2}\left(  h\right)  }%
{nr_{n}^{2}v^{2}\left(  h\right)  }+\frac{v^{2}\left(  h\right)  }{r_{n}^{2}%
}.
\]

\end{lemma}

\textbf{Proof : }Since $\Gamma_{n,K}^{\dag}$ is a positive operator, its
square root exists and :%
\begin{align*}
\left\langle \Gamma_{n,K}^{\dag}\overline{Z}_{K,n},\overline{Z}_{K,n}%
\right\rangle  &  =\left\Vert \left(  \Gamma_{n,K}^{\dag}\right)
^{1/2}\overline{Z}_{K,n}\right\Vert ^{2}\\
&  \leq C\left[  \left\Vert \left(  \Gamma_{n,K}^{\dag}\right)  ^{1/2}\left(
\overline{Z}_{K,n}-\mathbb{E}\overline{Z}_{K,n}\right)  \right\Vert
^{2}+\left\Vert \left(  \Gamma_{n,K}^{\dag}\right)  ^{1/2}\mathbb{E}%
\overline{Z}_{K,n}\right\Vert ^{2}\right]  .
\end{align*}
Then :%
\begin{align*}
&  \left\langle \Gamma_{n,K}^{\dag}\overline{Z}_{K,n},\overline{Z}%
_{K,n}\right\rangle ^{2}\\
&  \leq C\left[  \left\Vert \left(  \Gamma_{n,K}^{\dag}\right)  ^{1/2}\left(
\overline{Z}_{K,n}-\mathbb{E}\overline{Z}_{K,n}\right)  \right\Vert
^{4}+\left\Vert \left(  \Gamma_{n,K}^{\dag}\right)  ^{1/2}\mathbb{E}%
\overline{Z}_{K,n}\right\Vert ^{4}\right] \\
&  \leq C\left\Vert \Gamma_{n,K}^{\dag}\right\Vert _{\infty}^{2}\left(
\left\Vert \overline{Z}_{K,n}-\mathbb{E}\overline{Z}_{K,n}\right\Vert
^{4}+\left\Vert \mathbb{E}\overline{Z}_{K,n}\right\Vert ^{4}\right)  .
\end{align*}
From Lemma \ref{L2} and Lemma \ref{L3} we get :%
\[
\mathbb{E}\left\langle \Gamma_{n,K}^{\dag}\overline{Z}_{K,n},\overline
{Z}_{K,n}\right\rangle ^{2}\leq C\frac{h^{4}F^{2}\left(  h\right)  }%
{nr_{n}^{2}v^{2}\left(  h\right)  }+\frac{v^{4}\left(  h\right)  }{r_{n}%
^{2}v^{2}\left(  h\right)  }.
\]

\subsection{Derivation of the main results}

We start with a short and simple intermezzo about optimization in Hilbert spaces.

\textbf{Proof of Proposition \ref{progtheorik} :}

Consider the program :%
\[
\min_{a\in\mathbb{R},\varphi\in H}\mathbb{E}\left[  \left(  y-a-\left\langle
\varphi,X-x_{0}\right\rangle \right)  ^{2}K\left(  \dfrac{\left\Vert
X-x_{0}\right\Vert }{h}\right)  \right]  .
\]
Simple computations lead to :%
\begin{align*}
\mathcal{E}\left(  a,\varphi\right)   &  =\mathbb{E}\left[  \left(
y-a-\left\langle \varphi,X-x_{0}\right\rangle \right)  ^{2}K\left(
\dfrac{\left\Vert X-x_{0}\right\Vert }{h}\right)  \right] \\
&  =C+a^{2}\mathbb{E}K+\left\langle \Gamma_{K}\varphi,\varphi\right\rangle
-2a\mathbb{E}\left(  yK\right)  -2\left\langle \mathbb{E}\left(  yZK\right)
,\varphi\right\rangle +2a\left\langle \mathbb{E}\left(  ZK\right)
,\varphi\right\rangle .
\end{align*}
Obviously $\mathcal{E}\left(  a,\varphi\right)  $ is positive strictly convex
and
\[
\lim_{a,\left\Vert \varphi\right\Vert \rightarrow+\infty}\mathcal{E}\left(
a,\varphi\right)  =+\infty
\]
hence $\mathcal{E}\left(  a,\varphi\right)  $ has a single minimum (see
Rockafellar (1996) for further information about the minimization of convex
functions). It is also differentiable for all $\left(  a,\varphi\right)  $ in
$\mathbb{R}\times H.$ We compute its gradient~:%
\[
\nabla\mathcal{E}\left(  a,\varphi\right)  =\left(
\begin{array}
[c]{c}%
2a\mathbb{E}K-2\mathbb{E}\left(  yK\right)  +2\left\langle \mathbb{E}\left(
ZK\right)  ,\varphi\right\rangle \\
2\Gamma_{K}\varphi-2\mathbb{E}\left(  yZK\right)  +2a\mathbb{E}\left(
ZK\right)
\end{array}
\right)
\]
from which we get the solutions $\left(  a^{\ast},\varphi^{\ast}\right)  $ :%
\[
\left(
\begin{array}
[c]{c}%
a^{\ast}\mathbb{E}K+\left\langle \mathbb{E}\left(  ZK\right)  ,\varphi^{\ast
}\right\rangle =\mathbb{E}\left(  yK\right) \\
\Gamma_{K}\varphi^{\ast}=\mathbb{E}\left(  yZK\right)  -a^{\ast}%
\mathbb{E}\left(  ZK\right)
\end{array}
\right)  .
\]
We see from the second line that $\varphi^{\ast}$ is not uniquely defined if
$\Gamma_{K}$ is not one to one. Taking $\varphi^{\ast}=\Gamma_{K}^{-1}\left(
\mathbb{E}\left(  yZK\right)  -a^{\ast}\mathbb{E}\left(  ZK\right)  \right)  $
we get $\widehat{m}_{n}\left(  x_{0}\right)  $ as announced.

The forthcoming Lemma assesses that the random denominator of our estimate may
be replaced by a non-random one.

\begin{lemma}
\label{L7}When both $\frac{h^{4}}{nr_{n}^{2}v^{2}\left(  h\right)  }$and
$\frac{v\left(  h\right)  }{r_{n}F\left(  h\right)  }$ tend to zero, the
following holds :%
\[
\dfrac{\sum_{i=1}^{n}\omega_{i,n}}{nK\left(  1\right)  F\left(  h\right)
}-1\overset{L^{2}}{\rightarrow}0.
\]

\end{lemma}

\textbf{Proof :}%

\[
\sum_{i=1}^{n}\omega_{i,n}=\sum_{i=1}^{n}K_{i}-n\left\langle \Gamma
_{n,K}^{\dag}\overline{Z}_{K,n},\overline{Z}_{K,n}\right\rangle
\]
hence%
\[
\dfrac{\sum_{i=1}^{n}\omega_{i,n}}{nK\left(  1\right)  F\left(  h\right)
}-1=\frac{\sum_{i=1}^{n}K_{i}}{nK\left(  1\right)  F\left(  h\right)
}-1-\frac{\left\langle \Gamma_{n,K}^{\dag}\overline{Z}_{K,n},\overline
{Z}_{K,n}\right\rangle }{F\left(  h\right)  }.
\]
From Lemmas \ref{L5} and \ref{L6} we deduce that the announced Lemma \ref{L7} holds.

\subsubsection{Variance term}

We study first (see \ref{variance}) : $T_{v,n}=\dfrac{\sum_{i=1}^{n}\left(
y_{i}-m\left(  X_{i}\right)  \right)  \omega_{i,n}}{\sum_{i=1}^{n}\omega
_{i,n}}=\dfrac{\sum_{i=1}^{n}\varepsilon_{i}\omega_{i,n}}{\sum_{i=1}^{n}%
\omega_{i,n}}.$ It is plain that $\mathbb{E}T_{v,n}=0$. Denote $\widetilde
{T}_{v,n}=\dfrac{\sum_{i=1}^{n}\varepsilon_{i}\omega_{i,n}}{nK\left(
1\right)  F\left(  h\right)  }.$ We have :%
\[
T_{v,n}-\widetilde{T}_{v,n}=T_{v,n}\left(  \frac{F\left(  h\right)  -\frac
{1}{nK\left(  1\right)  }\sum_{i=1}^{n}\omega_{i,n}}{F\left(  h\right)
}\right)  .
\]
We begin with a Proposition. By Lemma \ref{L7} just above we know that
$T_{v,n}\sim\widetilde{T}_{v,n}$ in $L^{2}$ sense i.e.%
\[
\frac{\widetilde{T}_{v,n}}{T_{v,n}}\overset{L^{2}}{\rightarrow}1.
\]
.

\begin{proposition}
\label{v}We have :%
\[
\mathbb{E}T_{v,n}^{2}\leq C\frac{1}{nF^{2}\left(  h\right)  }\left(  F\left(
h\right)  +\frac{h^{2}F\left(  h\right)  }{nr_{n}v\left(  h\right)  }%
+\frac{v\left(  h\right)  }{r_{n}}\right)  .
\]
\bigskip
\end{proposition}

\textbf{Proof :} As announced above it suffices to prove the Proposition for
$\widetilde{T}_{v,n}$.%
\begin{align*}
\mathbb{E}\widetilde{T}_{v,n}^{2}  &  =\mathbb{E}\left(  \dfrac{\sum_{i=1}%
^{n}\varepsilon_{i}\omega_{i,n}}{nK\left(  1\right)  F\left(  h\right)
}\right)  ^{2}\\
&  =\frac{1}{n^{2}K^{2}\left(  1\right)  F^{2}\left(  h\right)  }%
\mathbb{E}\left\{  \mathbb{E}\left[  \left(  \sum_{i=1}^{n}\varepsilon
_{i}\omega_{i,n}\right)  ^{2}|X_{1},...,X_{n}\right]  \right\} \\
&  =\frac{1}{n^{2}K^{2}\left(  1\right)  F^{2}\left(  h\right)  }%
\mathbb{E}\left[  \mathbb{E}\left(  \sum_{i=1}^{n}\varepsilon_{i}^{2}%
\omega_{i,n}^{2}|X_{1},...,X_{n}\right)  \right]
\end{align*}
since for $i\neq j$ $\mathbb{E}\left[  \left(  \varepsilon_{i}\omega
_{i,n}\varepsilon_{j}\omega_{j,n}\right)  |X_{1},...,X_{n}\right]
=\omega_{i,n}\omega_{j,n}\mathbb{E}\left[  \left(  \varepsilon_{i}%
\varepsilon_{j}\right)  |X_{1},...,X_{n}\right]  =0.$ Hence :%
\[
\mathbb{E}\widetilde{T}_{v,n}^{2}=\frac{1}{n^{2}K^{2}\left(  1\right)
F^{2}\left(  h\right)  }\sigma_{\varepsilon}^{2}\mathbb{E}\left(  \sum
_{i=1}^{n}\omega_{i,n}^{2}\right)  =\frac{\sigma_{\varepsilon}^{2}%
\mathbb{E}\left(  \omega_{1,n}^{2}\right)  }{nK^{2}\left(  1\right)
F^{2}\left(  h\right)  }.
\]
By Lemma \ref{L4},%
\[
\mathbb{E}\left(  \omega_{1,n}^{2}\right)  \leq C\left(  F\left(  h\right)
+\frac{h^{2}F\left(  h\right)  }{nr_{n}v\left(  h\right)  }+\frac{v\left(
h\right)  }{r_{n}}\right)
\]
from which we deduce the Proposition.

Now we turn to the bias term.

\subsubsection{Bias term}

Remember that we have to deal with :%
\[
T_{b,n}=\dfrac{\sum_{i=1}^{n}\left(  m\left(  X_{i}\right)  -m\left(
x_{0}\right)  \right)  \omega_{i,n}}{\sum_{i=1}^{n}\omega_{i,n}}.
\]

Copying what was done above with $T_{v,n},$ we know that we can focus on :%
\[
\widetilde{T}_{b,n}=\dfrac{\sum_{i=1}^{n}\left(  m\left(  X_{i}\right)
-m\left(  x_{0}\right)  \right)  \omega_{i,n}}{nK\left(  1\right)  F\left(
h\right)  }%
\]
via Lemma \ref{L7}. For each $i$ there exists $c_{i}\in B\left(
x_{0},h\right)  $ such that :%
\begin{align*}
&  m\left(  X_{i}\right)  -m\left(  x_{0}\right) \\
&  =\left\langle m^{\prime}\left(  x_{0}\right)  ,Z_{i}\right\rangle +\frac
{1}{2}\left\langle m^{\prime\prime}\left(  c_{i}\right)  \left(  Z_{i}\right)
,Z_{i}\right\rangle .
\end{align*}
with $Z_{i}=X_{i}-x_{0}$. We deal with the first and second order derivatives
separatedly : $\widetilde{T}_{b,n}=\widetilde{T}_{b,n,1}+\widetilde{T}%
_{b,n,2}$ with
\begin{align*}
\widetilde{T}_{b,n,1}  &  =\frac{\sum_{i=1}^{n}\left\langle m^{\prime}\left(
x_{0}\right)  ,Z_{i}\right\rangle \omega_{i,n}}{nK\left(  1\right)  F\left(
h\right)  },\\
\widetilde{T}_{b,n,2}  &  =\frac{1}{2}\frac{\sum_{i=1}^{n}\left\langle
m^{\prime\prime}\left(  c_{i}\right)  \left(  Z_{i}\right)  ,Z_{i}%
\right\rangle \omega_{i,n}}{nK\left(  1\right)  F\left(  h\right)  }.
\end{align*}

\begin{proposition}
\label{b1}We have :%
\[
\mathbb{E}\widetilde{T}_{b,n,1}^{2}\leq C\frac{h^{2}}{nF\left(  h\right)
}+C\frac{v^{2}\left(  h\right)  }{F^{2}\left(  h\right)  }.
\]

\end{proposition}

\textbf{Proof of the Proposition :}

We first see that :%
\begin{align*}
\sum_{i=1}^{n}\left\langle m^{\prime}\left(  x_{0}\right)  ,X_{i}%
-x_{0}\right\rangle \omega_{i,n}  &  =\sum_{i=1}^{n}\left\langle m^{\prime
}\left(  x_{0}\right)  ,Z_{i}\right\rangle K_{i}\left(  1-\left\langle
Z_{i},\Gamma_{n,K}^{\dag}\overline{Z}_{K,n}\right\rangle \right) \\
&  =\sum_{i=1}^{n}\left\langle m^{\prime}\left(  x_{0}\right)  ,Z_{i}%
\right\rangle K_{i}\\
&  -\sum_{i=1}^{n}\left\langle m^{\prime}\left(  x_{0}\right)  ,Z_{i}%
\right\rangle K_{i}\left\langle Z_{i},\Gamma_{n,K}^{\dag}\overline{Z}%
_{K,n}\right\rangle \\
&  =n\left\langle m^{\prime}\left(  x_{0}\right)  ,\overline{Z}_{K,n}%
\right\rangle -n\left\langle \Gamma_{n,K}m^{\prime}\left(  x_{0}\right)
,\Gamma_{n,K}^{\dag}\overline{Z}_{K,n}\right\rangle \\
&  =n\left\langle m^{\prime}\left(  x_{0}\right)  ,\left(  I-\Gamma
_{n,K}\Gamma_{n,K}^{\dag}\right)  \overline{Z}_{K,n}\right\rangle
\end{align*}
and%
\[
\widetilde{T}_{b,n,1}=\frac{\left\langle m^{\prime}\left(  x_{0}\right)
,\left(  I-\Gamma_{n,K}\Gamma_{n,K}^{\dag}\right)  \left(  \overline{Z}%
_{K,n}\right)  \right\rangle }{K\left(  1\right)  F\left(  h\right)  }.
\]

Then we split into two terms :%
\begin{align*}
\left\langle m^{\prime}\left(  x_{0}\right)  ,\left(  I-\Gamma_{n,K}%
\Gamma_{n,K}^{\dag}\right)  \left(  \overline{Z}_{K,n}\right)  \right\rangle
&  =\left\langle \left(  I-\Gamma_{n,K}\Gamma_{n,K}^{\dag}\right)  m^{\prime
}\left(  x_{0}\right)  ,\left(  \overline{Z}_{K,n}-\mathbb{E}\overline
{Z}_{K,n}\right)  \right\rangle \\
&  +\left\langle \left(  I-\Gamma_{n,K}\Gamma_{n,K}^{\dag}\right)  m^{\prime
}\left(  x_{0}\right)  ,\mathbb{E}\overline{Z}_{K,n}\right\rangle .
\end{align*}
The $L^{2}$ norm of the first is bounded by $Ch\sqrt{F\left(  h\right)  /n}$
(see Lemma \ref{L3}) and the $L^{2}$ norm of the second is bounded by
$Cv\left(  h\right)  $ (see Lemma \ref{L2}). This finishes the proof of
Proposition \ref{b1}.

We turn to $\widetilde{T}_{b,n,2}$ and cut it into two parts :%
\begin{align*}
\widetilde{T}_{b,n,2}  &  =\frac{1}{2}\frac{\sum_{i=1}^{n}\left\langle
m^{\prime\prime}\left(  c_{i}\right)  \left(  Z_{i}\right)  ,Z_{i}%
\right\rangle K_{i}}{nK\left(  1\right)  F\left(  h\right)  }\\
&  -\frac{1}{2}\frac{\sum_{i=1}^{n}\left\langle m^{\prime\prime}\left(
c_{i}\right)  \left(  Z_{i}\right)  ,Z_{i}\right\rangle K_{i}\left\langle
Z_{i},\Gamma_{n,K}^{\dag}\overline{Z}_{K,n}\right\rangle }{nK\left(  1\right)
F\left(  h\right)  }\\
&  =R_{bn1}+R_{bn2}.
\end{align*}
The two forthcoming Propositions aim at giving a bound for the mean square
norm of $R_{bn1}$ and $R_{bn2}$.

\begin{proposition}
\label{b2}We get :%
\[
\mathbb{E}R_{bn1}^{2}\leq C\left(  \frac{h^{4}}{nF\left(  h\right)  }%
+h^{4}\right)  .
\]

\end{proposition}

\textbf{Proof of the Proposition :}

It is plain to see that for all $i$ and when Assumption $\mathbf{A}_{5}$ holds
:%
\[
0\leq\left\langle m^{\prime\prime}\left(  c_{i}\right)  \left(  Z_{i}\right)
,Z_{i}\right\rangle K_{i}\leq\left(  \mathit{\sup_{x\in\mathcal{V}\left(
x_{0}\right)  }}\left\Vert m^{\prime\prime}\left(  x\right)  \right\Vert
_{\infty}\right)  \left\Vert Z_{i}\right\Vert ^{2}K_{i}%
\]
hence that :%
\[
0\leq R_{bn1}\leq\frac{C}{2}\frac{\sum_{i=1}^{n}\left\Vert Z_{i}\right\Vert
^{2}K_{i}}{nK\left(  1\right)  F\left(  h\right)  }%
\]

It follows that :%
\[
0\leq R_{bn1}^{2}\leq C\frac{\left(  \sum_{i=1}^{n}\left\Vert Z_{i}\right\Vert
^{2}K_{i}\right)  ^{2}}{n^{2}F^{2}\left(  h\right)  }.
\]
Then :%
\begin{align*}
0  &  \leq\mathbb{E}R_{bn1}^{2}\leq\frac{C}{F^{2}\left(  h\right)  }%
\mathbb{E}\left(  \frac{1}{n}\sum_{i=1}^{n}K_{i}\left\Vert Z_{i}\right\Vert
^{2}\right)  ^{2}\\
&  =\frac{C}{F^{2}\left(  h\right)  }\left[  \frac{1}{n}\mathbb{E}\left(
K_{i}^{2}\left\Vert Z_{i}\right\Vert ^{4}\right)  +\frac{1}{n^{2}}\sum_{1\leq
i\neq j\leq n}\mathbb{E}\left(  K_{i}\left\Vert Z_{i}\right\Vert ^{2}%
K_{j}\left\Vert Z_{j}\right\Vert ^{2}\right)  \right] \\
&  \leq\frac{C}{F^{2}\left(  h\right)  }\left[  \frac{1}{n}\mathbb{E}K_{i}%
^{2}\left\Vert Z_{i}\right\Vert ^{4}+\left(  \mathbb{E}K_{i}\left\Vert
Z_{i}\right\Vert ^{2}\right)  ^{2}\right] \\
&  \leq\frac{C}{F^{2}\left(  h\right)  }\left[  \frac{h^{4}F\left(  h\right)
}{n}+h^{4}F^{2}\left(  h\right)  \right] \\
&  =C\left[  \frac{h^{4}}{nF\left(  h\right)  }+h^{4}\right]  .
\end{align*}

We turn to $R_{bn2}$.

\begin{proposition}
\label{b3}We have :%
\[
\mathbb{E}R_{bn2}^{2}\leq C\frac{h^{6}}{r_{n}^{2}}.
\]

\end{proposition}

Dealing with $R_{bn2}$ is a bit more complicated. We have :%
\[
-2R_{bn2}=\frac{1}{K\left(  1\right)  F\left(  h\right)  }\left\langle
\Gamma_{n,K}^{\dag}\overline{Z}_{K,n},\frac{1}{n}\sum_{i=1}^{n}\left\langle
m^{\prime\prime}\left(  c_{i}\right)  \left(  Z_{i}\right)  ,Z_{i}%
\right\rangle K_{i}Z_{i}\right\rangle .
\]
The next operation consists in replacing $\overline{Z}_{K,n}$ by its
expectation. Like above in the proof of Proposition \ref{b1} as well as in the
proof of Lemmas \ref{L4} and \ref{L6} we can add and subtract $\mathbb{E}ZK$
from $\overline{Z}_{K,n}$. Once again we decide not to go through details here
for the sake of shortness and clarity. Finally since the remaining involving
$\overline{Z}_{K,n}-\mathbb{E}\overline{Z}_{K,n}$ tends to zero quicker in
mean square, we can focus on :%
\begin{equation}
4R_{bn2}^{2}\leq\frac{C}{F^{2}\left(  h\right)  }\left\Vert \Gamma_{n,K}%
^{\dag}\right\Vert _{\infty}^{2}\left\Vert \mathbb{E}KZ\right\Vert
^{2}\left\Vert \frac{1}{n}\sum_{i=1}^{n}\left\langle m^{\prime\prime}\left(
c_{i}\right)  \left(  Z_{i}\right)  ,Z_{i}\right\rangle K_{i}Z_{i}\right\Vert
^{2}. \label{zeend}%
\end{equation}
At last we have to deal with :%
\[
\mathbb{E}\left\Vert \frac{1}{n}\sum_{i=1}^{n}\left\langle m^{\prime\prime
}\left(  c_{i}\right)  \left(  Z_{i}\right)  ,Z_{i}\right\rangle K_{i}%
Z_{i}\right\Vert ^{2}.
\]
Easy computations give :%
\begin{align}
&  \left\Vert \frac{1}{n}\sum_{i=1}^{n}\left\langle m^{\prime\prime}\left(
c_{i}\right)  \left(  Z_{i}\right)  ,Z_{i}\right\rangle K_{i}Z_{i}\right\Vert
^{2}\nonumber\\
&  =\frac{1}{n^{2}}\sum_{i=1}^{n}\left\langle m^{\prime\prime}\left(
c_{i}\right)  \left(  Z_{i}\right)  ,Z_{i}\right\rangle ^{2}K_{i}%
^{2}\left\Vert Z_{i}\right\Vert ^{2}\nonumber\\
&  +\frac{2}{n^{2}}\sum_{i<j}\left\langle m^{\prime\prime}\left(
c_{i}\right)  \left(  Z_{i}\right)  ,Z_{i}\right\rangle \left\langle
m^{\prime\prime}\left(  c_{j}\right)  \left(  Z_{j}\right)  ,Z_{j}%
\right\rangle \left\langle K_{i}Z_{i},K_{j}Z_{j}\right\rangle . \label{cfjp}%
\end{align}
We take expectations now and apply assumption $\mathbf{A}_{5}$ to the first
sum :%
\begin{align*}
&  \frac{1}{n^{2}}\mathbb{E}\sum_{i=1}^{n}\left\langle m^{\prime\prime}\left(
c_{i}\right)  \left(  Z_{i}\right)  ,Z_{i}\right\rangle ^{2}K_{i}%
^{2}\left\Vert Z_{i}\right\Vert ^{2}\\
&  \leq C\frac{1}{n^{2}}\mathbb{E}\sum_{i=1}^{n}\left\Vert Z_{i}\right\Vert
^{4}K_{i}^{2}\left\Vert Z_{i}\right\Vert ^{2}\\
&  =\frac{C}{n}\mathbb{E}\left(  K_{i}^{2}\left\Vert Z_{i}\right\Vert
^{6}\right)  \leq\frac{C}{n}h^{6}F\left(  h\right)  .
\end{align*}
Since $h^{6}F\left(  h\right)  /n$ tends to zero at a rate much quicker than
the next term we do not let it appear in the Proposition.

We fix $i$ and $j$ in (\ref{cfjp}) and take expectation :%
\begin{align*}
&  \mathbb{E}\left\langle m^{\prime\prime}\left(  c_{i}\right)  Z_{i}%
,Z_{i}\right\rangle \left\langle m^{\prime\prime}\left(  c_{j}\right)
Z_{j},Z_{j}\right\rangle \left\langle K_{i}Z_{i},K_{j}Z_{j}\right\rangle \\
&  =\left\langle \mathbb{E}\left[  \left\langle m^{\prime\prime}\left(
c_{i}\right)  Z_{i},Z_{i}\right\rangle K_{i}Z_{i}\right]  ,\mathbb{E}\left[
\left\langle m^{\prime\prime}\left(  c_{j}\right)  Z_{j},Z_{j}\right\rangle
K_{j}Z_{j}\right]  \right\rangle \\
&  =\left\Vert \mathbb{E}\left[  \left\langle m^{\prime\prime}\left(
c\right)  Z,Z\right\rangle KZ\right]  \right\Vert ^{2}.
\end{align*}
By assumption $A_{5}$ we get :%
\begin{align*}
&  \left\vert \mathbb{E}\left\langle m^{\prime\prime}\left(  c_{i}\right)
Z_{i},Z_{i}\right\rangle \left\langle m^{\prime\prime}\left(  c_{j}\right)
\left(  Z_{j}\right)  ,Z_{j}\right\rangle \left\langle K_{i}Z_{i},K_{j}%
Z_{j}\right\rangle \right\vert \\
&  \leq\left(  \mathbb{E}\left\Vert \left\langle m^{\prime\prime}\left(
c\right)  Z,Z\right\rangle KZ\right\Vert \right)  ^{2}\\
&  \leq C\left[  \mathbb{E}\left(  K\left\Vert Z\right\Vert ^{3}\right)
\right]  ^{2}\\
&  \leq Ch^{6}F^{2}\left(  h\right)  .
\end{align*}
Finally with (\ref{zeend}) at hand we have :%
\begin{align*}
\mathbb{E}R_{bn2}^{2}  &  \leq\frac{C}{F^{2}\left(  h\right)  }\frac{v^{2}%
}{v^{2}r_{n}^{2}}\left(  \frac{C}{n}h^{6}F\left(  h\right)  +Ch^{6}%
F^{2}\left(  h\right)  .\right) \\
&  \leq C\frac{h^{6}}{r_{n}^{2}}%
\end{align*}
since $nF\left(  h\right)  \rightarrow+\infty$.

At last we finish with the proof of the main Theorem which is considerably
alleviated by all that was done above.\bigskip

\textbf{Proof of Theorem \ref{main}, Proposition \ref{MA} and Proposition
\ref{opt} :}

The proof of the Theorem stems from display (\ref{moon}), Propositions
\ref{v}, \ref{b1}, \ref{b2} and \ref{b3}. Collecting these previous results we
have :%
\begin{align*}
\mathbb{E}\left(  \widehat{m}_{n}\left(  x_{0}\right)  -m\left(  x_{0}\right)
\right)  ^{2}  &  \leq C\frac{1}{nF^{2}\left(  h\right)  }\left(  F\left(
h\right)  +\frac{h^{2}F\left(  h\right)  }{nr_{n}v\left(  h\right)  }%
+\frac{v\left(  h\right)  }{r_{n}}\right) \\
&  +C\left[  \frac{h^{6}}{r_{n}^{2}}+h^{4}+\frac{h^{2}}{nF\left(  h\right)
}+\frac{v^{2}\left(  h\right)  }{F^{2}\left(  h\right)  }\right]  .
\end{align*}
First from%
\[
v\left(  h\right)  \leq h^{2}F\left(  h\right)  ,
\]
we see that the first line above will be an $O\left(  1/\left(  nF\left(
h\right)  \right)  \right)  $ whenever $h^{2}/r_{n}$ and $h^{2}/\left(
nr_{n}v\left(  h\right)  \right)  $ are bounded. We turn to the second line.
The term is at least $h^{2}/\left(  nF\left(  h\right)  \right)  $ may be
removed because it can be neglected with repect to the variance term. In order
to reach an $O\left(  h^{4}\right)  $ for the bias we have to bound
$h^{2}/r_{n}^{2}$ and $1/\left(  h^{2}nF\left(  h\right)  \right)  $.

At last summing up all what was done above comes down to taking $r_{n}\asymp
h,$ and $n\cdot\min\left\{  v\left(  h\right)  /h,h^{2}F\left(  h\right)
\right\}  \geq C>0$.

Following the results of Mas (2007b) this last inequality comes down, when
$\rho$ is regularly varying at $0$ with positive index :%
\[
nF\left(  h\right)  \cdot\min\left\{  \rho\left(  h\right)  ,h^{2}\right\}
\geq C>0
\]
And Theorem \ref{main} is proved.

\end{document}